\documentclass[10pt]{amsart}
\usepackage{amsfonts}
\usepackage{indentfirst}
\usepackage{graphicx}
\usepackage{amssymb}
\usepackage[colorlinks,citecolor=red,urlcolor=blue,bookmarks=false,hypertexnames=true]{hyperref}
\usepackage{color,xcolor}
\usepackage{bm}
\usepackage{multirow}
\usepackage{multicol}
\usepackage{amssymb,amsmath,amsthm,mathrsfs}%
\usepackage{epsfig}
\usepackage{epstopdf}
\usepackage{subfigure,caption}
\usepackage{abstract}

\allowdisplaybreaks

\newtheorem{theorem}{Theorem}[section]

\theoremstyle{definition}

\theoremstyle{remark}
\newtheorem{remark}[theorem]{Remark}
\newtheorem{lemma}{Lemma}[section]

\numberwithin{equation}{section}%
\numberwithin{table}{section}%
\numberwithin{figure}{section}

\def\3bar{{|\hspace{-.02in}|\hspace{-.02in}|}}

\def\tc{\text{curl}}

\def\RR{\mathbb{R}}

\def\NN{\mathbb{N}}

\def\RR{\mathbb{R}}

\def\NN{\mathbb{N}}
\def\curl{\text{curl}}

\def\QQ{\mathbb{Q}}

\def\tr{\mathsf{T}}
\def\bs{\boldsymbol}

\def\curl{\operatorname{curl}}

\newcommand{\OO}[1]{\overset{\scriptscriptstyle \circ}{#1}{}}

\graphicspath{{fig/fig1/}{fig/fig2/}{fig/}{figure/}}

\begin{document}
	\title	[H($\tc^2$)-conforming quadrilateral spectral element method]{H($\tc^2$)-conforming  quadrilateral spectral element method  for  quad-curl  problems}
	
	  \author[
	L. Wang,  W. Shan,   H. Li,  \& Z. Zhang
	]{Lixiu Wang${}^{1}$,  \;\;
		Weikun Shan${}^{2}$, \;\;  Huiyuan Li${}^{3}$, \;\;   and \;\; Zhimin Zhang${}^{4}$
		}
	
	\thanks{${}^{1}$Beijing Computational Science Research Center, Beijing 100193, China. Email: {\tt lxwang@csrc.ac.cn}.
\\
\indent ${}^{2}$School of Science, Henan University of Technology, Zhengzhou, Henan, 450001, China.  Email: {\tt shanweikun66@126.com}. The research of this author is supported   in part by
the National Natural Science Foundation of China (NSFC 11801147) and
 the  Fundamental Research Funds for the Henan Provincial Colleges and Universities in Henan University of Technology (2017QNJH20).
\\
\indent ${}^{3}$State Key Laboratory of Computer Science/Laboratory of Parallel Computing,  Institute of Software, Chinese Academy of Sciences, Beijing 100190, China. Email: {\tt huiyuan@iscas.ac.cn}. The research of this author is  supported in part by the National Natural Science Foundation of China (NSFC 11871145 and NSFC 11971016) and National Key R\&D Program of China (No. 2018YFB0204404).
\\
\indent ${}^{4}$Beijing Computational Science Research Center, Beijing 100193, China;  and Department of Mathematics, Wayne State University, MI 48202, USA. Email: {\tt zmzhang@csrc.ac.cn; ag7761@wayne.edu}. The research of this author is supported
in part by the National Natural Science Foundation of China (NSFC  11871092, NSFC 11926356, and NSAF U1930402).
}
	

%



\maketitle
\begin{abstract}
In this paper, we propose an $H(\tc^2)$-conforming  quadrilateral spectral element method to solve  quad-curl problems.  Starting with generalized Jacobi polynomials, we first introduce quasi-orthogonal polynomial systems for vector fields  over rectangles.
$H(\tc^2)$-conforming  elements over arbitrary  convex  quadrilaterals  are  then  constructed explicitly in a hierarchical pattern using the contravariant transform together with the bilinear mapping from the reference square onto each quadrilateral.  It is astonishing that
both the simplest   rectangular  and  quadrilateral spectral elements possess  only 8 degrees of freedom on each physical element.  In the sequel, we propose our $H(\tc^2)$-conforming quadrilateral spectral element approximation  based on
the mixed weak formulation  to solve the quad-curl equation and its eigenvalue problem.  Numerical results show the effectiveness and  efficiency of our method.
\end{abstract}

\section{Introduction}
Quad-curl problems, including the Maxwell's transmission eigenvalue problem (MTEP) \cite{Sun2016Finite, Cakoni2017A} and  the  resistive magnetohydrodynamic (MHD) \cite{Qingguo2012A, Zheng2011A}, have received increasing  attention in recent years. The transmission eigenvalue problem plays an important role  in  qualitative approaches for  the inverse scattering theory for inhomogeneous media, while MHD models have widely used in thermonuclear fusion, plasma physics, geophysics, and astrophysics \cite{Qingguo2012A}.
It is  meaningful and  urgent  to design  highly efficient and accurate  numerical methods for quad-curl problems.
 By the way,
singularities of the quad-curl problem were analyzed  in \cite{Nicaise2018Singularities, Zhang2018M2NA}.

In contrast to second-order curl problems, limited work has been done on numerical methods for quad-curl problems.  Initially,  numerical methods with various nonconformity/mix techniques,  such as  nonconforming finite element  methods \cite{Zheng2011A}, discontinuous Galerkin methods \cite{Qingguo2012A}, weak Galerkin methods \cite{Sun2019Awg}, mixed finite element methods \cite{Sun2016A,Sun_2013, Monk2012Finite, Wang2019A,Zhang2017Regular, Zhang2018M2NA,Sun2018Multigrid}, and the Hodge decomposition method \cite{Brenner2017Hodge,Brenner2019Multigrid}, were proposed to solve  quad-curl problems as well as  their related eigenvalue problems.
Indeed, $H(\tc^2)$-conforming methods were unavailable  for quad-curl problems
until recently.  In \cite{Zhang2019H}, $H(\tc^2)$-conforming finite elements  were first proposed over  parallelograms and  triangles with  their convergence analysis being carried out both theoretically and numerically.  Although  incomplete polynomials are adopted to reduce the  number of basis functions,  this conforming method still has 24 degrees of freedom (DOFs) on each parallelogram element.  In addition, even for the lowest order $H(\tc^2)$-conforming element, the construction of the  Lagrange type basis functions are very complicated, which prevents the implement for higher-order elements. In another recent work \cite{WangToappear}, in order to solve quad-curl eigenvalue problems, a family of $H(\tc^2)$-conforming finite elements  over triangles are constructed using  complete polynomials of total degree $\ge 4$
 with at least $30$ DOFs on each element. More recently,
  by introducing  continuous and discrete de Rham complexes with high order Sobolev spaces,  Hu et al. discovered in \cite{Hu2020simple}  that the simplest rectangular finite  element  possess only 8 DOFs.
 Until now, there is no $H(\tc^2)$-conforming finite elements  designed for general quadrilateral meshes and no systematic way to construct $H(\tc^2)$-conforming elements of arbitrarily high orders.
This paper  is then motivated by the desire of  $H(\tc^2)$-conforming elements  of an arbitrarily high order over arbitrary  convex quadrilaterals. 


As one of the most  important high order methods, the spectral element method was first introduced by Patera \cite{Patera1984A}. In analogy to $p$- and
$hp$-finite element methods, spectral element methods inherit the high-order  convergence of
the traditional spectral methods, while preserve the
flexibility of the low-order finite element methods \cite{Li2019C1}.
There is abundant literature  addressing spectral/$hp$ element approximations  for  second-order electromagnetic equations (see  \cite{Belgacem1999Spectral, Cohen2002Higher, Sabine2006,Liu2010A, Na2015Mixed, wangzhangzhang2018, Cai2013} and the reference therein), which validates the superiority of spectral/$hp$ element   methods over low order methods.
However, no efforts have been reported in literature on $H(\tc^2)$-conforming spectral element methods up to now.
Indeed,
more stringent  continuity requirements should be imposed on $H(\tc^2)$-conforming spectral elements than $H(\tc)$-conforming elements, which hinder the progress on  the construction  of the $H(\tc^2)$-conforming basis functions,   especially those over quadrilaterals.



The aim of the current paper is to construct hierarchical $H(\tc^2)$-conforming basis functions  on general quadrangulated meshes, and then to propose   an efficient quadrilateral spectral element approximation for solving quad-curl problems directly. Similar to $C^1$-conforming basis functions \cite{Li2019C1}, $H(\tc^2)$-conforming basis functions can be divided into vertex modes, edge modes, and interior modes.  The interior modes are  constructed such that their tangential components and  their curls   are zeros along every edge. The edge modes involve eight one-dimensional trace functions constituted of function values and curls on four edges. For each edge basis function, all trace functions but one vanish identically. Inspired by \cite{Zhang2019H}, the vertex modes on a quadrilateral adopt 4 DOFs determined  by their curls at four
vertices. 

Based on the de Rham complex,
we first introduce  quasi-orthogonal polynomial systems for vector fields on the reference square  with the help of  generalized Jacobi polynomials  of indexes $(-1,-1)$ and $(-2,-2)$. Their tangent and curl traces along four edges are then explored.
 Accompanying the bilinear mapping from the reference square to each quadrilateral element, we introduce a  contravariant  transformation between vector fields on the reference  square and those on the physical quadrilateral, which preserve  their tangent components along edges up to some constants  and their curls up to the Jacobian of  the bilinear mapping.
This characteristic gives us a quick and easy construction of the interior modes and  the tangent edge modes on a physical element from those  vectorial polynomial basis on the reference square whose curls are zero along four edges.  The curl edge modes  on the physical element are then derived by simply multiplying the corresponding vectorial polynomial basis  on the reference square.
While vertex modes are also  technically  set up by using  those vectorial polynomial basis  belonging to $Q_{2,3}\times Q_{3,2}$  on the reference square.
We note that basis functions on a quadrilateral will reduce to those on a rectangle whenever the  quadrilateral falls into a rectangle.
It is  then not surprise that the simplest element has 8 DOFs  on a quadrilateral just  as the simplest  element on a rectangle.
With the help of our $H(\curl^2)$-conforming spectral elements, we  finally propose an efficient and direct  quadrilateral spectral element method to solve the quad-curl problems.

The rest of the paper is organized as follows. In Section 2 we list some function spaces and notations. $H(\tc^2)$-conforming spectral elements over  arbitrary quadrilaterals are defined in Section 3.  Special  attention is paid to those  over parallelograms or rectangles.
Section 4  is devoted to the technical  derivation of  our $H(\tc^2)$-conforming spectral elements.   In Section 5,  we propose the $H(\tc^2)$-conforming spectral element method to solve the quad-curl problem with the divergence constraint
on the  basis of its mixed weak formulation. In Section 6, numerical examples are presented  to  verify the  correctness and efficiency of our method. Some concluding remarks  are finally  given in  Section 7.

\section{Preliminaries}
\subsection{Notations}
Denote by $\NN_0$,  $\NN$ and $\RR$  the collections of non-negative integers, positive integers and real numbers, respectively.
Let $\Omega\subset\mathbb{R}^2$ be a convex Lipschitz domain,  and $\bm n$ be the unit outward normal vector to $\partial \Omega$. We adopt standard notations for Sobolev spaces such as $H^m(\Omega)$ or $H_0^m(\Omega)$ with the norm $\left\|\cdot\right\|_{m,\Omega}$ and the semi-norm $\left|\cdot\right|_{m,\Omega}$. If $m=0$,  the space $H^0 (\Omega)$ coincides with $ L^2(\Omega)$ equipped with the norm $\|\cdot\|_{\Omega}$. We shall drop the subscript $\Omega$ whenever no confusion would arise. We use  ${\bm L}^2(\Omega)$ to denote the vector-valued Sobolev spaces  $L^2(\Omega)^2$.
Further we denote by  $\QQ_{m,n}(\Omega)$  the bivariate polynomial space
of separate degrees at most $m$ and $n$.
%

Let ${\bm u}=(u_1, u_2)^{\tr}$ and ${\bm w}=(w_1, w_2)^{\tr}$, where the superscript ${\tr}$ denotes the transpose.
Then ${\bm u} \times {\bm w} = u_1 w_2 - u_2 w_1$ and $\nabla \times {\bm u} = \frac{\partial u_2}{\partial x} - \frac{\partial u_1}{\partial y}$.
For a scalar function $v$, $\nabla \times v = (\frac{\partial v}{\partial y}, - \frac{\partial v}{\partial x})^{\tr}$.
We denote $(\nabla\times)^2\bm u=\nabla\times\nabla\times\bm u$ and  define
\begin{align*}
H(\text{curl};\Omega)&:=\{\bm u \in {\bm L}^2(\Omega):\; \nabla \times \bm u \in L^2(\Omega)\},\\
H(\text{curl}^2;\Omega)&:=\{\bm u \in {\bm L}^2(\Omega):\; \nabla \times \bm u \in L^2(\Omega),\;(\nabla \times)^2 \bm u \in \bm L^2(\Omega)\},
\end{align*}
whose norms are defined by
\[\left\|\bm u\right\|^2_{H(\tc^s;\Omega)}=\sum_{i=0}^s\|(\nabla\times)^i\bm u\|_\Omega^2\]
 with $s=1,\,2.$
The spaces $H_0^s(\text{curl};\Omega)(s=1,\,2)$ are defined as follows:
\begin{align*}
&H_0(\text{curl};\Omega):=\{\bm u \in H(\text{curl};\Omega):\;{\bm n}\times\bm u=0\; \text{on}\ \partial \Omega\},\\
&H_0(\text{curl}^2;\Omega):=\{\bm u \in H(\text{curl}^2;\Omega):\;{\bm n}\times\bm u=0\; \text{and}\; \nabla\times \bm u=0\;\; \text{on}\ \partial \Omega\}.
\end{align*}


 \subsection{Generalized Jacobi polynomials}
We introduce some generalized Jacobi polynomials which play an important role in constructing  the $H(\tc^2)$-conforming elements.
For any $\alpha,\beta>-1,\ n\in \mathbb{N}_0$, denote $J_n^{\alpha,\beta}(\zeta)$ as the $n$-th classic Jacobi polynomial with respect to the weight function $(1-\zeta)^\alpha(1+\zeta)^\beta$ on $[-1,1]$. Various generalization have been introduced to allow $\alpha$ and/or $\beta$ being negative integers \cite{Guo2009Generalized,Guo2010Composite,Szeg1939Orthogonal}. In this paper, we use the following generalized Jacobi polynomials:
 \begin{align}\label{Jacobidefinition}
{K}^{-1,-1}_{n}(\zeta)=\begin{cases}
\dfrac{1-\zeta}{2}, &n=0,\\[0.5em]
\dfrac{1+\zeta}{2}, &n=1,\\[0.5em]
\dfrac{\zeta^2-1}{4}J_{n-2}^{1,1}(\zeta),&n\geq 2.
\end{cases}
\qquad
{K}^{-2,-2}_{n}(\zeta)=\begin{cases}
\dfrac{(1-\zeta)^2(2+\zeta)}{4},&n=0,\\[0.5em]
\dfrac{(1-\zeta)^2(1+\zeta)}{4}, &n=1,\\[0.5em]
\dfrac{(1+\zeta)^2(2-\zeta)}{4}, &n=2,\\[0.5em]
\dfrac{(1+\zeta)^2(\zeta-1)}{4}, &n=3,\\[0.5em]
\left(\dfrac{\zeta^2-1}{4}\right)^2J_{n-4}^{2,2}(\zeta), &n\geq 4.
\end{cases}
\end{align}

It is readily checked that $\{ K^{-1,-1}_n: n\ge 2\}$ and $\{ K^{-2,-2}_n: n\ge 4\}$
coincide, up to a constant, with the generalized Jacobi polynomials defined in \cite{Guo2009Generalized, Shen2009A},
while
$\{ K^{-1,-1}_n: 0\le n\le 1\}$ and $\{ K^{-2,-2}_n: 0\le n\le 3\}$
are exactly  Lagrange and Hermite   interpolating basis functions on $[-1,1]$,
respectively.  More precisely,  for $n\in \mathbb{N}_0$,
\begin{align}
\label{Hermite1}
&  K^{-1,-1}_n(-1) =  \delta_{0, n} , \qquad\quad  K^{-1,-1}_n(1) =  \delta_{1, n},
\\
\label{Hermite2}
&\partial_{\zeta}^l  K^{-2,-2}_n(-1) =  \delta_{l, n} , \qquad
\partial_{\zeta}^l  K^{-2,-2}_n(1) =  \delta_{l+2, n} , \qquad   0\le l\le 1.
\end{align}
Besides, it holds that
\begin{align}
\label{Hermite4}
&K^{-1,-1}_n{}'(\zeta)  = \dfrac{n-1}{2}K^{0,0}_{n-1}(\zeta),\quad  n\geq2,\\
\label{Hermite5}
& K^{-2,-2}_n{}'(\zeta)  = \dfrac{n-3}{2}K^{-1,-1}_{n-1}(\zeta),\quad n\geq4,\\
\label{Hermite5-01}
& K^{-2,-2}_0{}'(\zeta)  = \dfrac{3(-1+\zeta)(\zeta+1)}{4},\quad K^{-2,-2}_1{}'(\zeta)=\dfrac{(-1+\zeta)(3 \zeta+1)}{4},\\
& K^{-2,-2}_2{}'(\zeta)  = \frac{3(1-\zeta)(\zeta+1)}{4},\quad K^{-2,-2}_3{}'(\zeta)=\frac{(\zeta+1)(-1+3 \zeta)}{4},\nonumber\\
\label{Hermite6}
&K^{-1,-1}_n(\zeta)=\dfrac{n-1}{2(2n-1)}\left(J^{0,0}_{n}(\zeta)-J^{0,0}_{n-2}(\zeta)\right),\quad  n\geq2.
\end{align}
Hereafter, we use the notation ${}'$ to denote  $\partial_\zeta$ when there is no confusion.
\subsection{Continuity  across $H(\text{curl}^2)$-conforming elements}
We  introduce the following lemma which shows the request of the  continuity across the cells' edges  of the $H(\text{curl}^2)$-conforming elements. It  has great significance for constructing the basis functions.

\begin{lemma}\label{prob2}\cite{Zhang2019H}
	Let $\text{K}_1$ and $\text{K}_2$ be two non-overlapping Lipschitz domains having a common edge $\Lambda$ such that $\overline{\text{K}_1}\cap\overline{\text{K}_2} = \Lambda$. Assume that  ${\bm u}_i=\bm u|_{K_i}\in H(\mathrm{curl}^2;\text{K}_i), i=1,2$, and  define
	\begin{equation*}
	\displaystyle{\bm u}=
	\begin{cases}
	\bm u_1,  & \text{in}\ \text{K}_1,\\
	\bm u_2, & \text{in}\ \text{K}_2.
	\end{cases}
	\end{equation*}
	Then $\bm u_1 \times \bm n_1 = -\bm u_2 \times \bm n_2$\ and $\nabla\times \bm u_1=\nabla\times \bm u_2$ on $\Lambda$ implies that $\bm u \in H(\mathrm{curl}^2;\text{K}_1 \cup \text{K}_2 \cup \Lambda)$, where $\bm n_i$ ($i=1,2$) is the unit outward normal vector to $K_i$ on $\Lambda$.
\end{lemma}

\section{$H(\text{curl}^2)$-conforming quadrilateral spectral elements }\label{quadrilateralelem}

In this section, we present our main theory about  the $H(\tc^2)$-conforming basis functions on
an arbitrary  quadrilateral. All basis functions are derived from the generalized Jacobi polynomials, and they are divided into vertex modes, edge modes and interior modes.  The tangential components  and curls of an interior mode are identically zero on all edges. The edge modes
are further divided into two groups: the first group (function edge modes) are  constructed  such that their tangential components  have magnitudes only on  one edge while
their curls are enforced zero on all edges;  the tangential components  of the second group (curl edge modes) vanish identically on all
edges while their curls have magnitudes  only on one edge. The vertex modes are devised such that their tangential components vanish identically on all
edges while their curls have magnitudes only on adjacent edges.


\subsection{Mapping from the reference square onto quadrilateral}
\begin{figure}[htbp]
\hfill
\begin{minipage}[htbp]{0.35\textwidth}
\centering
\includegraphics[width=1\textwidth]{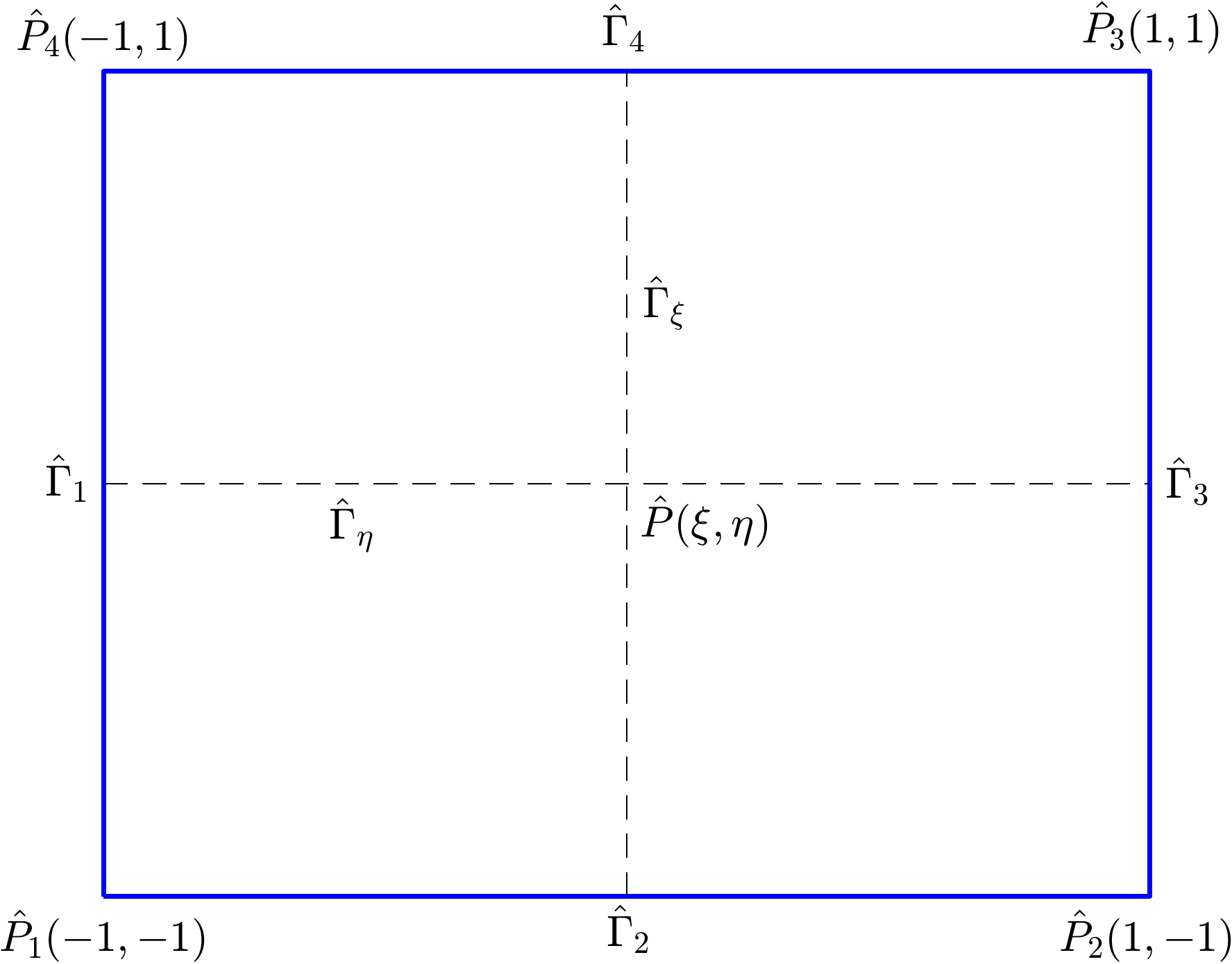}
\end{minipage}
\hfill
\begin{minipage}[htbp]{0.35\textwidth}
\centering
\includegraphics[width=1\textwidth]{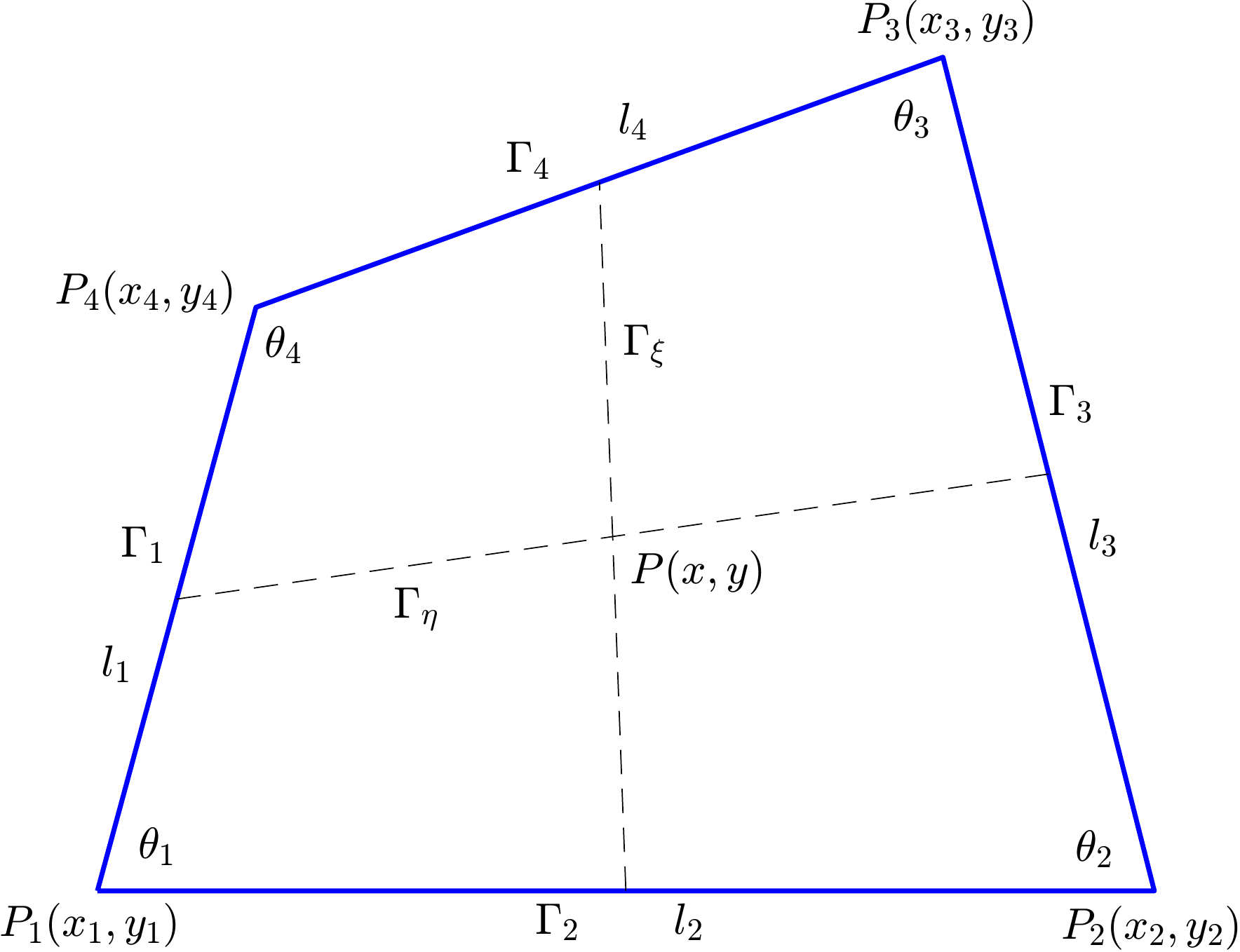}
\end{minipage}
\hspace*{\fill}
\caption{Reference square $\hat{K}$(left) and arbitrary convex quadrilateral(right).}
\label{domain-1}
\end{figure}

 Let $\hat{K}=(-1,1)^2$ be the reference square with the vertices $\hat{P}_i$, and the edges $\hat{\Gamma}_i, 1\leq i\leq 4$. Let $K$ be an arbitrary convex quadrilateral with the vertices $P_i(x_i, y_i)$, the edges $\Gamma_i$, and the inner angles $\theta_i, 1\leq i\leq 4$; see Figure \ref{domain-1}. The side length of $\Gamma_i$ is denote by $l_i$.  We  emphasize that the tangential directions of $\hat {\bm\tau}$ and  $\bm \tau$ are anti-clockwise. Follow the line in \cite{Li2019C1}, we define the one-to-one mapping $\Phi_K: \hat K\mapsto K$ such that

\begin{align}\label{tran2}
\begin{split}
\begin{pmatrix}
    x \\
    y  \end{pmatrix}
=\Phi_K\begin{pmatrix}
    \hat x \\
    \hat y  \end{pmatrix}:&=\sigma_1(\hat x,\hat y)\begin{pmatrix}
    x_1 \\
    y_1
  \end{pmatrix}
+\sigma_2(\hat x,\hat y)\begin{pmatrix}
    x_2 \\
    y_2
     \end{pmatrix}
    \\&
+\sigma_3(\hat x,\hat y)
\begin{pmatrix}
    x_3 \\
    y_3
     \end{pmatrix}
+\sigma_4(\hat x,\hat y)
\begin{pmatrix}
    x_4 \\
    y_4
     \end{pmatrix}, \qquad (\hat x, \hat y)\in \bar{\hat K}, \, ( x,  y)\in  \bar K,
\end{split}
\end{align}
where \begin{align*}
&\sigma_1(\hat x,\hat y)=\frac{(1-\hat{x})(1-\hat{y})}{4},\ \sigma_2(\hat x,\hat y)=\frac{(1+\hat{x})(1-\hat{y})}{4},\\
 &\sigma_3(\hat x,\hat y)=\frac{(1+\hat{x})(1+\hat{y})}{4},\ \sigma_4(\hat x,\hat y)=\frac{(1-\hat{x})(1+\hat{y})}{4}.
\end{align*}
For simplicity, we write $x_{ij}=x_i-x_j,\ y_{ij}=y_i-y_j$ hereafter. Given a scalar function $\phi$ defined on $K$, we associate it with $ {\hat{\phi}} :=\phi \circ \Phi_{K}$ on $\hat K$.
It is easy to see from the chain rule that
\begin{align}\label{chain_rule}
\partial_{\hat x}\hat \phi =B_{11}(\hat y) \partial_{x}\phi +B_{21}(\hat y) \partial_{y}\phi,\quad \partial_{\hat y}\hat \phi =B_{12}(\hat x) \partial_{x}\phi+B_{22}(\hat x) \partial_{y}\phi,
\end{align}
where
\begin{align*}
&B_{11}=B_{11}(\hat y) :=\frac{x_{21}}{2} \frac{1-\hat y}{2}+\frac{x_{34}}{2} \frac{1+\hat y}{2}, \quad B_{12}=B_{12}(\hat x) :=\frac{x_{41}}{2} \frac{1-\hat x}{2}+\frac{x_{32}}{2} \frac{1+\hat x}{2},\\ &B_{21}=B_{21}(\hat y) :=\frac{y_{21}}{2} \frac{1-\hat y}{2}+\frac{y_{34}}{2} \frac{1+\hat y}{2}, \quad B_{22}=B_{22}(\hat x) :=\frac{y_{41}}{2} \frac{1-\hat x}{2}+\frac{y_{32}}{2} \frac{1+\hat x}{2}.
\end{align*}

 The Jacobian matrix of the transformation $\Phi_K$ with respect to the reference coordinates and its determinant are denoted by
 \begin{align}
 \label{FK}
 B_K(\hat{x},\hat{y}) =
\begin{pmatrix}
      \dfrac{\partial x}{\partial \hat{x}} & \dfrac{\partial x}{\partial \hat{y}} \\[0.7em]
      \dfrac{\partial y}{\partial \hat{x}} & \dfrac{\partial y}{\partial \hat{y}}
   \end{pmatrix}=\begin{pmatrix}
      B_{11} & B_{12} \\[0.7em]
      B_{21} & B_{22}
   \end{pmatrix},
\qquad J_K(\hat{x},\hat{y})=\det(B_K(\hat{x},\hat{y})), \qquad (\hat{x},\hat{y})\in {\hat {K}}.
\end{align}
More clearly,
\begin{equation}
\begin{aligned} J_K(\hat x, \hat y) :&=\frac{x_{21} y_{41}-y_{21} x_{41}}{4} \sigma_{1}(\hat x, \hat y)+\frac{x_{32} y_{12}-y_{32} x_{12}}{4} \sigma_{2}(\hat x, \hat y) \\ &+\frac{x_{43} y_{23}-y_{43} x_{23}}{4} \sigma_{3}(\hat x, \hat y)+\frac{x_{14} y_{34}-y_{14} x_{34}}{4} \sigma_{4}(\hat x, \hat y) \\ =&\,\frac{l_{2} l_{1} \sin \theta_{1}}{4} \sigma_{1}(\hat x, \hat y)+\frac{l_{3} l_{2} \sin \theta_{2}}{4} \sigma_{2}(\hat x, \hat y)+\frac{l_{4} l_{3} \sin \theta_{3}}{4} \sigma_{3}(\hat x, \hat y)+\frac{l_{1} l_{4} \sin \theta_{4}}{4} \sigma_{4}(\hat x, \hat y). \end{aligned}
\end{equation}

Based on the bilinear mapping \eqref{tran2},  we introduce  the contravariant transformation from $H(\tc^2;\hat K)$ to $H(\tc^2; K)$,
\begin{align}
 \bs{p}  =  B_K^{-\tr}  \hat {\bs{p}}  \circ \Phi_K^{-1}, \qquad \hat {\bs{p}} \in H(\tc^2;\hat K).\label{trans-01}
 \end{align}
It can be checked from \eqref{chain_rule} that
\begin{align}
& \nabla \times \bs {p} = J^{-1}_K  \hat \nabla  \times \hat{\bs{p}} \circ \Phi_K^{-1},\label{trans-02}
 \\
&(\nabla \times)^2 \bs {p} =\dfrac{B_K}{J_K^2}\left((\hat\nabla \times)^2 \hat{\bm p} -\dfrac{\hat\nabla \times \hat{\bm p}}{J_K}
\begin{pmatrix}
      \dfrac{x_{12}+x_{34}}{4}B_{22}-\dfrac{y_{12}+y_{34}}{4}B_{12} \\[0.7em]
      -\dfrac{y_{12}+y_{34}}{4}B_{11}+\dfrac{x_{12}+x_{34}}{4}B_{21}
   \end{pmatrix}\right)  \circ \Phi_K^{-1},\label{trans-03}\\
&  {\bs{p}} \cdot \bs{\tau}  = B_K^{-\tr}  \hat {\bs{p}} \cdot \frac{ B_K \hat {\bs \tau} }{\|B_K \hat {\bs \tau}\|} \circ \Phi_K^{-1}
= \frac{  \hat {\bs{p} } \cdot \hat {\bs \tau} }{\|B_K \hat {\bs \tau}\|} \circ \Phi_K^{-1}.\label{trans-04}
 \end{align}
We note that  $\| B_K \hat {\bs{\tau}} \|=\frac{l_i}{2}$ on ${\Gamma_i}$ under the bilinear mapping  \eqref{tran2}.
Throughout this paper, we always associate a vector field $\bm \phi\in H(\tc^2;K)$ with $\hat{\bm \phi}\in H(\tc^2;\hat K)$ by the contravariant transformation \eqref{trans-01}.
\begin{remark}
	
Note that $\Phi_K$ is reduced to affine mapping such that   $B_K$ is a constant matrix and $J_K$ is  constant whenever $K$ is a rectangle or a parallelogram. In particular, if  $K$ is a rectangle with $\Gamma_1$ and $\Gamma_3$ (resp. $\Gamma_2$ and $\Gamma_4$)  parallel  to  the vertical  (resp. horizontal) axis,
$B_K=\mathrm{Diag}\big(\frac{l_2}{2}, \frac{l_1}{2}\big)$ is a constant diagonal matrix.

\end{remark}
\subsection{$H(\tc^2)$-conforming quadrilateral spectral elements}

Now we are in a position to present hierarchical basis functions
of $H(\tc^2)$-conforming spectral elements on an unstructured quadrilateral mesh.
For simplicity, we shall write $s_i:=\sin \theta_i$, $i=1,2,3,4$.
We emphasize that our desired basis functions,  $\bs\psi_{m,n}, \bs\phi_{m,n}$, $m,n\ge 0$, on a quadrilateral  element $K$
are obtained from the functions, $\hat{\bs\psi}_{m,n}, \hat{\bs\phi}_{m,n}$, $m,n\ge 0$,  in the reference coordinates  through  the bilinear mapping \eqref{tran2} and the contravariant transformation \eqref{trans-01} such that
$\bs\psi_{m,n} = B_K^{-\tr}  \hat {\bs\psi}_{m,n} \circ \Phi_K^{-1}$, and $\bs\phi_{m,n} = B_K^{-\tr}  \hat {\bs\phi}_{m,n} \circ \Phi_K^{-1}$.

Details on the constructions of the $H(\tc^2)$-conforming quadrilateral elements will be postponed to the next section.

\begin{itemize}
\item Interior modes
\begin{align}
\label{Qinterior}
\begin{split}
&\hat {\bs \phi}_{m,n} =  \hat \nabla \big[K^{-1,-1}_m(\hat x) K^{-1,-1}_n(\hat y)\big],  \quad  m,n\ge 2,\\
& \hat {\bs\psi}_{m,n} =  \hat \nabla K^{-2,-2}_m(\hat x)  K^{-2,-2}_n(\hat y),   \quad \  \ m\in\{2,4,5,6,\cdots\}, n\ge 4,\\
& \hat {\bs\psi}_{m,2} = \hat K^{-2,-2}_m(\hat x)  \hat\nabla\hat  K^{-2,-2}_n(\hat y), \quad \  \  m\ge 4, n= 2,
\end{split}
\end{align}
such that
\begin{align*}
&{\bs \phi}_{m,n}\cdot {\bs \tau} \big|_{\Gamma} =  [0,0 , 0,0 ],
&&\nabla \times {\bs \phi}_{m,n} \big|_{\Gamma} =  [ 0  ,0,0,0],\\
&{\bs \psi}_{m,n}\cdot {\bs \tau} \big|_{\Gamma} =  [0,0 , 0,0 ],
&&\nabla \times {\bs \psi}_{m,n} \big|_{\Gamma} =  [ 0  ,0,0,0].
\end{align*}
Hereafter, we use tetrads for the trace on four edges, i.e., $\phi|_{\Gamma}=\left[\phi|_{\Gamma_1},\phi|_{\Gamma_2},\phi|_{\Gamma_3},\phi|_{\Gamma_4}\right]$.

\item  Function edge modes

\begin{itemize}
\item Function edge modes corresponding to $\Gamma_1$:
\begin{align}
\label{Qfunc1}
\begin{split}
&\hat {\bs \phi}_{0,n} =  \hat \nabla \big[K^{-1,-1}_0(\hat x) K^{-1,-1}_n(\hat y)\big], \quad n\ge 2,
\\
&\hat {\bs \phi}_{0,0}   =\left(\frac{\hat y(\hat y^2-1)\left(3 \hat x^{2}-5\right)}{32},-\frac{\hat x(\hat x^2-1)\left(3 \hat y^{2}-5\right)}{32}-\frac{\hat x-1}{4}\right)^\tr,
\end{split}
\end{align}
such that
\begin{align*}
&{\bs \phi}_{0,n}\cdot {\bs \tau} \big|_{\Gamma} =  \left[ -\frac{(n-1)J^{0,0}_{n-1}(\hat y) }{l_1},0 , 0,0 \right],
&&  \nabla \times {\bs \phi}_{0,n} \big|_{\Gamma} =  [ 0  ,0,0,0],\quad n\ge 2,\\
&{\bs \phi}_{0,0}\cdot {\bs \tau}  \big|_{\Gamma} = \left[ -\frac{1 }{l_1},0,0,0\right],
&&  \nabla \times {\bs \phi}_{0,0} \big|_{\Gamma} =  [ 0,0,0,0].
\end{align*}
\item  Function edge modes corresponding to $\Gamma_2$:
\begin{align}
\label{Qfunc2}
\begin{split}
&\hat {\bs \phi}_{m,0} =  \hat \nabla \big[K^{-1,-1}_m(\hat x) K^{-1,-1}_0(\hat y)\big], \quad   m\ge 2,
\\
&\hat {\bs \phi}_{1,0} =\left(-\frac{\hat y(\hat y^2-1)\left(3 \hat x^{2} -5\right)}{32}-\frac{\hat y-1}{4}, \frac{\hat x(\hat x^2-1)\left(3 \hat y^{2}-5\right)}{32}\right)^\tr,
\end{split}
\end{align}
such that
\begin{align*}
&{\bs \phi}_{m,0}\cdot {\bs \tau} \big|_{\Gamma} =  \left[ 0 , \frac{(m-1)J^{0,0}_{m-1}(\hat x) }{ l_2 },0,0 \right],
&& \nabla \times {\bs \phi}_{m,0} \big|_{\Gamma} =  [ 0  ,0,0,0],\quad m\ge 2,
\\
&{\bs \phi}_{1,0}\cdot {\bs \tau}  \big|_{\Gamma} = \left[0, \frac{1 }{l_2 },0,0\right],
&&  \nabla \times {\bs \phi}_{1,0} \big|_{\Gamma} =  [ 0,0,0,0].
\end{align*}
\item Function edge modes corresponding to $\Gamma_3$:
\begin{align}
\label{Qfunc3}
\begin{split}
&\hat {\bs \phi}_{1,n} =  \hat \nabla \big[K^{-1,-1}_1(\hat x) K^{-1,-1}_n(\hat y)\big], \quad   n\ge 2,
\\
&\hat {\bs \phi}_{1,1}= \left(-\frac{\hat y(\hat y^2-1)\left(3 \hat x^{2}-5\right)}{32}, \frac{\hat x(\hat x^2-1)\left(3 \hat y^{2} -5\right)}{32}+\frac{1+\hat x}{4}\right)^\tr,
\end{split}
\end{align}
such that
\begin{align*}
&{\bs \phi}_{1,n}\cdot {\bs \tau} \big|_{\Gamma} =  \left[ 0 , 0,\frac{(n-1)J^{0,0}_{n-1}(\hat y) }{ l_3},0 \right],
&&  \nabla \times {\bs \phi}_{1,n} \big|_{\Gamma} =  [ 0  ,0,0,0], \quad n\ge 2,\\
&{\bs \phi}_{1,1}\cdot {\bs \tau}  \big|_{\Gamma} = \left[0,0, \frac{1 }{ l_3 },0\right],
&& \nabla \times {\bs \phi}_{1,1} \big|_{\Gamma} =  [ 0,0,0,0].
\end{align*}

\item Function edge modes corresponding to $\Gamma_4$:
\begin{align}
\label{Qfunc4}
\begin{split}
&\hat {\bs \phi}_{m,1} =  \hat \nabla \big[K^{-1,-1}_m(\hat x) K^{-1,-1}_1(\hat y)\big], \quad   m\ge 2,
\\
&\hat {\bs \phi}_{0,1} =\left(\frac{\hat y(\hat y^2-1)\left(3 \hat x^{2} -5\right)}{32}+\frac{1+\hat y}{4},-\frac{\hat x(\hat x^2-1)\left(3 \hat y^{2}-5\right)}{32}\right)^\tr,
\end{split}
\end{align}
such that
\begin{align*}
&{\bs \phi}_{m,1}\cdot {\bs \tau} \big|_{\Gamma} =  \left[ 0  ,0,0,- \frac{(m-1)J^{0,0}_{m-1}(\hat x) }{ l_4} \right],
&&\nabla \times {\bs \phi}_{m,1} \big|_{\Gamma} =  [ 0  ,0,0,0],\quad m\ge2,\\
&{\bs \phi}_{0,1}\cdot {\bs \tau}  \big|_{\Gamma} = \left[0,0,0,- \frac{1 }{l_4}\right],
&&  \nabla \times {\bs \phi}_{0,1} \big|_{\Gamma} =  [ 0,0,0,0].
\end{align*}

\end{itemize}

By adding a multiple of the interior mode  $\hat{\bs \phi}_{3,3}$ to each basis functions,   one obtains
four  alternative function edge modes $\hat{\tilde{\bm\phi}}_{0,0}, \, \hat{\tilde{\bm\phi}}_{1,0},\,  \hat{\tilde{\bm\phi}}_{1,1}$
and $\hat{\tilde{\bm\phi}}_{0,1}$ with even  simple presentations:
 \begin{align*}
&\hat{\tilde{\bm\phi}}_{0,0}
=\hat{  \bm\phi }_{0,0}+\frac{\hat{\bs \phi}_{3,3}}{8} = \left(\frac{3\hat y(\hat y^2-1)(\hat x^2-1)}{16},   \frac{(\hat x+2)(\hat x-1)^2}{8}\right)^{\tr},
\\
&\hat{\tilde{\bm\phi}}_{1,0}
=\hat{  \bm\phi }_{1,0}+\frac{\hat{\bs \phi}_{3,3}}{8}= \left(\frac{(\hat y+2)(\hat y-1)^2}{8}, \frac{3x(x^2-1)(\hat y^2-1)}{16}\right)^{\tr},
\\
&\hat{\tilde{\bm\phi}}_{1,1}
=\hat{  \bm\phi }_{1,1}-\frac{\hat{\bs \phi}_{3,3}}{8}=\left(-\frac{3\hat y(\hat y^2-1)(\hat x^2-1)}{16},   -\frac{(\hat x-2)(\hat x+1)^2}{8}\right)^{\tr},
\\
&\hat{\tilde{\bm\phi}}_{0,1}
=\hat{  \bm\phi }_{0,1}-\frac{\hat{\bs \phi}_{3,3}}{8}= \left( -\frac{(\hat y-2)(\hat y+1)^2}{8}, -\frac{3\hat x(\hat x^2-1)(\hat y^2-1)}{16} \right)^{\tr}.
\end{align*}


\item Curl boundary modes

\begin{itemize}
\item Curl edge modes corresponding to $\Gamma_1$:
\begin{align}
&\bs{\hat  \psi}_{1,n} = J_K(\hat x, \hat y)  K^{-2,-2}_1(\hat x) \hat \nabla  K^{-2,-2}_n(\hat y),
\quad n\in\{2,4,5,6,\cdots\},\label{gamma1curl1}
\end{align}
such that
\begin{align*}
&{\bs \psi}_{1,n}\cdot {\bs \tau}  \big|_{\Gamma} = [0,0,0,0],
&& \nabla \times {\bs \psi}_{1,n} \big|_{\Gamma} = [ K^{-2,-2}_n{}'(\hat y),0,0,0].
\end{align*}

\item Curl edge modes corresponding to $\Gamma_2$:
\begin{align}
&\bs{\hat  \psi}_{m,1} = J_K(\hat x, \hat y) \hat \nabla K^{-2,-2}_m(\hat x) K^{-2,-2}_1(\hat y),
\quad m\in\{2,4,5,6,\cdots\},\label{gamma2curl1}
\end{align}
such that
\begin{align*}
&{\bs \psi}_{m,1}\cdot {\bs \tau}  \big|_{\Gamma} = [0,0,0,0],
&&\nabla  \times {\bs \psi}_{m,1} \big|_{\Gamma} =  [ 0,-K^{-2,-2}_m{}'(\hat x),0,0].
\end{align*}

\item Curl edge modes corresponding to $\Gamma_3$:
\begin{align}
&\bs{\hat  \psi}_{3,n} = J_K(\hat x, \hat y)  K^{-2,-2}_3(\hat x) \hat \nabla  K^{-2,-2}_n(\hat y),
\quad n\in\{2,4,5,6,\cdots\},\label{gamma3curl1}
\end{align}
such that
\begin{align*}
&{\bs \psi}_{3,n}\cdot {\bs \tau}  \big|_{\Gamma} = [0,0,0,0],
&&\nabla  \times {\bs \psi}_{3,n} \big|_{\Gamma} = [ 0,0,K^{-2,-2}_n{}'(\hat y),0].
\end{align*}

\item Curl edge modes corresponding to $\Gamma_4$:
\begin{align}
&\bs{\hat  \psi}_{m,3} = J_K(\hat x, \hat y) \hat \nabla K^{-2,-2}_m(\hat x) K^{-2,-2}_3(\hat y),
\quad m\in\{2,4,5,6,\cdots\},\label{gamma4curl1}
\end{align}
such that
\begin{align*}
&{\bs \psi}_{m,3}\cdot {\bs \tau}  \big|_{\Gamma} = [0,0,0,0],
&& \nabla \times {\bs \psi}_{m,3} \big|_{\Gamma} =    [ 0,0,0,-K^{-2,-2}_m{}'(\hat x)].
\end{align*}
\end{itemize}

Whenever $K$ is a rectangle or a parallelogram, $J_K(\hat x, \hat y)=\frac{|K|}{4}$, where $|K|$ stands for the area of $K$.


\item Vertex modes

\begin{itemize}
\item Vertex modes corresponding to $P_1$:
\begin{align}  \label{Qvert1}  
\begin{split}
\bm{\hat  \psi}_{0,0} =\Bigg(&(\hat y-1)^{2}(1+\hat y)(\hat x-1)\left(\frac{l_{2}\left(l_{1} s_{1}+2 l_{3} s_{2}\right) \hat x}{128}+\frac{l_{2}\left(3 l_{1} s_{1}+2 l_{3} s_{2}\right)}{128}\right),
\\
&-(1+\hat x)(\hat y-1)(\hat x-1)^{2}\left(\frac{l_{1}\left(l_{2} s_{1}+2 l_{4} s_{4}\right) \hat y}{128}+\frac{l_{1}\left(3 l_{2} s_{1}+2 l_{4} s_{4}\right)}{128}\right)\Bigg)^\tr,
\end{split}
\end{align}
such that
\begin{align*}
& {\bs \psi}_{0,0}\cdot {\bs \tau}  \big|_{\Gamma} = [0,0,0,0],
&& \nabla \times {\bs \psi}_{0,0} \big|_{\Gamma} = \left[-\frac{\hat{y}}{2}+\frac{1}{2},\frac{1}{2}-\frac{\hat{x}}{2}, 0,0\right].
\end{align*}

\item Vertex modes corresponding to $P_2$:
\begin{align}
\label{Qvert2}
\begin{split}
\bm{\hat  \psi}_{0,1} =\Bigg(&(\hat y-1)^{2}(1+\hat y)(\hat x+1)\left(\frac{l_{2}\left(2 l_{1} s_{1}+l_{3} s_{2}\right) \hat x}{128}-\frac{l_{2}\left(2 l_{1} s_{1}+3 l_{3} s_{2}\right)}{128}\right),\\
&(1-\hat x)(\hat y-1)(1+\hat x)^{2}\left(\frac{l_{3}\left(l_{2} s_{2}+2 l_{4} s_{3}\right) \hat y}{128}+\frac{l_{3}\left(3 l_{2} s_{2}+2 l_{4} s_{3}\right)}{128}\right)\Bigg)^\tr,
\end{split}
\end{align}
such that
\begin{align*}
&  {\bs \psi}_{0,1}\cdot {\bs \tau}  \big|_{\Gamma} = [0,0,0,0],
&& \nabla \times {\bs \psi}_{0,1} \big|_{\Gamma} = \left[0,\frac{1}{2}+\frac{\hat{x}}{2}, -\frac{\hat{y}}{2}+\frac{1}{2}, 0\right].
\end{align*}

\item Vertex modes corresponding to $P_3$:
\begin{align}
\label{Qvert3}
\begin{split}
\bm{\hat  \psi}_{1,1} =\Bigg(&(\hat y+1)^{2}(1-\hat y)(\hat x+1)\left(-\frac{l_{4}\left(2 l_{1} s_{4}+l_{3} s_{3}\right) \hat x}{128}+\frac{l_{4}\left(2 l_{1} s_{4}+3 l_{3} s_{3}\right)}{128}\right),\\
&(1-\hat x)(\hat y+1)(\hat x+1)^{2}\left(\frac{l_{3}\left(2 l_{2} s_{2}+l_{4} s_{3}\right) \hat y}{128}-\frac{l_{3}\left(2 l_{2} s_{2}+3 l_{4} s_{3}\right)}{128}\right)\Bigg)^\tr,
\end{split}
\end{align}
such that
\begin{align*}
& {\bs \psi}_{1,1}\cdot {\bs \tau}  \big|_{\Gamma} = [0,0,0,0],
&& \nabla \times {\bs \psi}_{1,1} \big|_{\Gamma} = \left[0,0, \frac{1}{2}+\frac{\hat{y}}{2}, \frac{1}{2}+\frac{\hat{x}}{2}\right].
\end{align*}

\item Vertex modes corresponding to $P_4$:
\begin{align}
\label{Qvert4}
\begin{split}
\bm{\hat  \psi}_{1,0} =&
\Bigg((\hat y+1)^{2}(1-\hat y)(\hat x-1)\left(-\frac{l_{4}\left(l_{1} s_{4}+2 l_{3} s_{3}\right) \hat x}{128}-\frac{l_{4}\left(3 l_{1} s_{4}+2 l_{3} s_{3}\right)}{128}\right), \\
&(1+\hat x)(\hat y+1)(\hat x-1)^{2}\left(-\frac{l_{1}\left(2 l_{2} s_{1}+l_{4} s_{4}\right) \hat y}{128}+\frac{l_{1}\left(2 l_{2} s_{1}+3 l_{4} s_{4}\right)}{128}\right)\Bigg)^\tr,
\end{split}
\end{align}
such that
\begin{align*}
& {\bs \psi}_{1,0}\cdot {\bs \tau}  \big|_{\Gamma} = [0,0,0,0],
&&  \nabla \times {\bs \psi}_{1,0} \big|_{\Gamma} = \left[\frac{1}{2}+\frac{\hat{y}}{2}, 0,0, \frac{1}{2}-\frac{\hat{x}}{2}\right].
\end{align*}

\end{itemize}

Whenever $K$ is a rectangle or parallelogram, $\bm{\hat  \psi}_{0,0},\bm{\hat  \psi}_{0,1},\bm{\hat  \psi}_{1,1},\bm{\hat  \psi}_{1,0}$ are reduced to
\begin{align*}
&\bm{\hat  \psi}_{0,0} =\left(\frac{|K| (\hat y-1)^{2}(\hat y+1)(\hat x-1)\left(3 \hat x+5\right)}{128},
-\frac{|K|  (\hat x+1)(\hat y-1)(\hat x-1)^{2}\left(3 \hat y+5\right)}{128}\right)^\tr,
\\
&\bm{\hat  \psi}_{0,1} =\left(\frac{|K|  (\hat y-1)^{2}(\hat y+1)(\hat x+1)\left(3\hat x-5\right)}{128}, -\frac{|K| (\hat x-1)(\hat y-1)(\hat x+1)^{2}\left(3\hat y+5\right)}{128}\right)^\tr,
\\
&\bm{\hat  \psi}_{1,1} =\Bigg(\frac{|K| (\hat y+1)^{2}(\hat y-1)(\hat x+1)\left(3\hat x-5\right)}{128},
-\frac{|K| (\hat x-1)(\hat y+1)(\hat x+1)^{2}\left(3\hat y-5\right)}{128}\Bigg)^\tr,
\\
&\bm{\hat  \psi}_{1,0} =\Bigg(\frac{|K| (\hat y+1)^{2}(\hat y-1)(\hat x-1)\left(3\hat x+5\right)}{128}, -\frac{|K| (\hat x+1)(\hat y+1)(\hat x-1)^{2}\left(3\hat y-5\right)}{128}\Bigg)^\tr,
\end{align*}
respectively.
\end{itemize}

\section{Constructions of the $H(\tc^2)$-conforming quadrilateral spectral elements}\label{construct}
In this section, we show the derivation process of $H(\tc^2)$-conforming basis functions on a quadrilateral element $K$
from basis functions for vector fields on the reference square $\hat K$  through the bilinear mapping  \eqref{tran2} and the contravariant  transformation \eqref{trans-01}.

 We  first resort to the following de Rham  complex \cite{Arnold2018Finite, Hu2020simple}
 for an enlightenment of the construction process,
  \begin{equation}\label{derham}
  \begin{aligned}
  &\mathbb{R}&\overset{id}\longrightarrow &H^1(\Omega)&\overset{\nabla}\longrightarrow &H(\tc^2;\Omega)&\overset{\nabla\times}\longrightarrow
  &\quad H^1(\Omega)&\overset{0}\longrightarrow &\{0\}\\
  &&&\quad \cup&&\quad \quad \cup&&\quad \quad \cup&\\
  &\mathbb{R}&\overset{id}\longrightarrow &\quad S_h&\overset{\nabla}\longrightarrow &\quad \quad W_h&\overset{\nabla\times}\longrightarrow &\quad \quad  U_h&\overset{0}\longrightarrow &\{0\},\\
  \end{aligned}
  \end{equation}
 where $S_h=\left\{ u_h\in H^1(\Omega): u_h\big|_{K} \circ\Phi^{-1}_K  \text{ is certain bivariate polynomial for all } K\in \mathcal{T}_h\right\}$, and $W_h, U_h$ are certain conforming  element spaces over the partition $\mathcal{T}_h$ of $\Omega$.
%
The main property of the de Rham complex lies in  that the range of each operator coincides with the kernel
of the following operator
consecutive operators. In particular,
 \begin{align}
 \label{curlkernel}
 \begin{aligned}
&H(\tc^2;\Omega)\hspace*{-1em} && \supseteq \left\{ \bs u\in H(\tc^2;\Omega):  \nabla \times \bs u =0\right\}\hspace*{-1em}  &&= \nabla H^1(\Omega) \hspace*{-1em} &&\supseteq \nabla H^2(\Omega),
\\
&W_h \hspace*{-1em} && \supseteq \left\{ \bs u\in W_h:  \nabla \times \bs u =0\right\}  &&= \nabla S_h &&\supseteq \nabla  \left[S_h\cap  H^2(\Omega)\right] .
 \end{aligned}
 \end{align}

\subsection{Vectorial basis functions on the reference square}
It is known that $ K^{-1,-1}_m(\hat x) K^{-1,-1}_n(\hat y)$,  $  m,n\geq 0$ are basis functions on the reference square $\hat K$
for  constructing $C^1$-conforming spectral elements using the bilinear mapping \eqref{tran2} over quadrilateral partitions \cite{Li2019C1}.
 According to \eqref{curlkernel} and \eqref{chain_rule},
 $\bs{B}_K^{-\tr} \hat \nabla\left[K^{-1,-1}_m(\hat x) K^{-1,-1}_n(\hat y) \right]\circ \Phi_K^{-1} ,
 \  m,n\geq 0,$
 are local basis candidates on $K$ for the discrete kernel space $\left\{ \bs u\in V_h:  \nabla \times \bs u =0\right\}$.
 For this reason, we would  like to choose
 \begin{align}
 \label{bases-01}
 \hat{\bm p}_{m,n}(\hat x,\hat y):=\hat \nabla \big[K^{-1,-1}_m(\hat x) K^{-1,-1}_n(\hat y)\big],\quad  (m,n)\in \NN_0^2\setminus\{(0,0)\}
 ,
 \end{align}
 as a part of vectorial basis functions on the reference square.

 Inspired by the success of the work in \cite{Sabine2006},
  we  recommend
%
   \begin{align}
   \label{bases-02}
  &\hat{\bm q}_{m,n}(\hat x,\hat y):=\hat \nabla K^{-2,-2}_m(\hat x) K^{-2,-2}_n(\hat y)\   \text{ or } \
   \hat{\bm q}_{m,n}^\ast(\hat x,\hat y):=K^{-2,-2}_m(\hat x) \hat \nabla K^{-2,-2}_n(\hat y), \quad  (m,n)\in \NN^2,
  \end{align}
for the other part of   vectorial basis functions on the reference square.

It is easy to see that  $\hat{\bm q}_{0,n}+\hat{\bm q}_{2,n}=\hat{\bm q}^\ast_{m,0}+\hat{\bm q}^\ast_{m,2}=0$
for any $m,n\in \NN_0$,
and
\begin{align*}
&\hat{\bm q}_{m,n}+\hat{\bm q}_{m,n}^\ast =\hat \nabla [\hat K^{-2,-2}_m(\hat x)  K^{-2,-2}_n(\hat y)],
\\
&\hat{\bm q}_{m,0}+\hat{\bm q}_{m,2} =\hat \nabla [\hat K^{-2,-2}_m(\hat x)],  &&
\hat{\bm q}^\ast_{0,n}+\hat{\bm q}^\ast_{2,n} =\hat \nabla [\hat K^{-2,-2}_n(\hat y)],
\end{align*}
are all finite linear combination of $\hat{\bs p}_{i,j}$, $(i,j)\in \NN_0^2\setminus\{(0,0)\}$ for any  $m,n\in \NN_0$.
Thus two parts of vectorial functions $\hat{\bs p}_{m,n}$, $(m,n)\in \NN_0^2\setminus\{(0,0)\}$,
 and $\hat{\bs q}_{m,n}/\hat{\bs q}^{\ast}_{m,n}$, $(m,n)\in \NN^2$, form a basis in $[L^2(\hat K)]^2$.

Now we are ready to present the trace properties of the  vectorial basis functions under
the contravariant transformation \eqref{trans-01}.  The following three lemmas 
can be deduced from the definition of the general Jacobi polynomials \eqref{Jacobidefinition}, and their properties \eqref{Hermite1}-\eqref{Hermite5-01} immediately.

\begin{lemma}\label{lemma-01-4}
Let $\bm{p}_{{m,n}}=B_K^{-\tr}\hat {\bm{p}}_{{m,n}}\circ \Phi_K^{-1}$, and it holds that
\begin{align}\label{lem-01-4-curl}
\nabla \times \bm{p}_{{m,n}} =&\, 0,
\end{align}
\begin{align}
\label{lem-01-4-tan}
\begin{split}
\bm{p}_{{m,n}} \cdot \bs{\tau} \big|_{\Gamma}
=  & \Bigg[ -\frac{  2K^{-1,-1}_m(-1) K^{-1,-1}_n{}'(\hat y) }{ l_1 } ,   \frac{  2K^{-1,-1}_m{}'(\hat x) K^{-1,-1}_n(-1)   }{ l_2 },\\
 &\quad \frac{  2K^{-1,-1}_m(1) K^{-1,-1}_n{}'(\hat y)   }{ l_3 },
  - \frac{2K^{-1,-1}_m{}'(\hat x) K^{-1,-1}_n(1)  }{ l_4 }\Bigg]
   \\
   = & \begin{cases}
     [0,0,0,0], & m\ge 2,\, n\ge 2,\\
    [ -\frac{(n-1)J^{0,0}_{n-1}(\hat y) }{ l_1 }  ,0,0,0], & m=0,\, n\ge 2,\\
    [ 0  ,0,\frac{(n-1)J^{0,0}_{n-1}(\hat y) }{ l_3 },0], & m=1,\, n\ge 2,\\
    [ 0,\frac{(m-1)J^{0,0}_{m-1}(\hat x) }{ l_2 }  ,0,0], & m\ge 2,\, n=0,\\
    [ 0  ,0,0,- \frac{(m-1)J^{0,0}_{m-1}(\hat x) }{ l_4 } ], & m\ge 2,\, n=1,\\
    [ \frac{1 }{ l_1 }  ,-\frac{1 }{ l_2 } , 0,0], & m=0,\, n=0, \\
    [ 0, \frac{1 }{ l_2 }  ,-\frac{1 }{ l_3 },0], & m=1,\, n=0, \\
    [ 0, 0 ,\frac{1 }{ l_3 },-\frac{1 }{ l_4 } ], & m=1,\, n=1,\\
    [ -\frac{1 }{ l_1 }, 0,0  ,\frac{1 }{ l_4 }], & m=0,\, n=1.
   \end{cases}
   \end{split}
\end{align}
\end{lemma}

\begin{lemma}\label{lemma-02-4}
Let ${\bm{q}}_{{m,n}}=B_K^{-\tr}\hat {\bm{q}}_{{m,n}}\circ \Phi_K^{-1}$, then it holds that
\begin{align}
\label{lem-02-4-tan}
\begin{split}
 {\bm{q}}_{{m,n}} \cdot \bm{\tau} \big|_{\Gamma}
= &\,  \left[0 ,   \frac{ 2 K^{-2,-2}_m{}'(\hat x) K^{-2,-2}_n(-1)   }{ l_2 },
0,
  - \frac{2K^{-2,-2}_m{}'(\hat x) K^{-2,-2}_n(1)  }{ l_4 }\right]
   \\
    = &\begin{cases}
     [0,0,0,0], & m\ge 1,\, n\in\{1,3,4,5,6,\cdots \} ,\\
    \left[0 ,   0,   0,   -\frac{2K^{-2,-2}_m{}'(\hat x)  }{ l_4 }\right], & m\ge 1,\, n=2,
   \end{cases}
\end{split}
 \end{align}
\begin{align}
\label{lem-02-4-curl}
\begin{split}
 \nabla \times {\bm{q}}_{{m,n}}|_{\Gamma} =&\, 
 -J^{-1}_K(\hat x,\hat y)   K^{-2,-2}_m{}'(\hat x)  K^{-2,-2}_n{}'(\hat y)\big|_{\Gamma}
 \\
 = &\begin{cases}
    [ 0,0,0,0], & m,n\in \{2,4,5,6,\cdots\},
    \\
     \left[-\frac{ K^{-2,-2}_n{}'(\hat y)}{J_K(-1, \hat y)},  0,0,0\right],  & m=1, \ n\in  \{2,4,5,6,\cdots\},
     \\
     \left[0,0,-\frac{K^{-2,-2}_n{}'(\hat y)}{J_K(1, \hat y)},  0\right],  & m=3, \ n\in  \{2,4,5,6,\cdots\},
     \\
     \left[0,-\frac{K^{-2,-2}_m{}'(\hat x)}{J_K(\hat x, -1)},  0,0\right],  & n=1,\  m\in  \{2,4,5,6,\cdots\},
     \\
     \left[0,0,0,-\frac{K^{-2,-2}_m{}'(\hat x)}{J_K(\hat x, 1)} \right],  & n=3, \ m\in  \{2,4,5,6,\cdots\},
     \\
     \left[-  \frac{K^{-2,-2}_1{}'(\hat y)}{J_K(-1,\hat y)}, -  \frac{K^{-2,-2}_1{}'(\hat x)}{J_K(\hat x,-1)} ,0, 0\right],  & m=1, \ n=1,
     \\
     \left[0, - \frac{K^{-2,-2}_3{}'(\hat x)}{J_K(\hat x,-1) } , -  \frac{K^{-2,-2}_1{}'(\hat y)}{J_K(1,\hat y)}, 0\right],  & m=3, \ n=1,
     \\
     \left[0 ,0, - \frac{K^{-2,-2}_3{}'(\hat y)}{J_K(1,\hat y) }, -\frac{K^{-2,-2}_3{}'(\hat x)}{J_K(\hat x,1)  }\right],  & m=3, \ n=3,
     \\
     \left[- \frac{K^{-2,-2}_3{}'(\hat y)}{J_K(-1,\hat y) }, 0, 0, - \frac{K^{-2,-2}_1{}'(\hat x)}{J_K(\hat x,1) } \right],  & m=1, \ n=3.
 \end{cases}
 \end{split}
 \end{align}
 \end{lemma}

\begin{lemma}\label{lemma-03-4}
Let\  ${\bm{q}}_{{m,n}}^\ast=B_K^{-\tr}\hat {\bm{q}}_{{m,n}}^\ast\circ \Phi_K^{-1}$, then we have
\begin{align}
\label{lem-03-4-tan}
\begin{split}
 {\bs{q}_{{m,n}}^\ast} \cdot \bs{\tau} \big|_{\Gamma}
= &\,  \left[-\frac{ 2 K^{-2,-2}_m{}(-1) K^{-2,-2}_n{}'(\hat y)   }{ l_1 },
0,
   \frac{2K^{-2,-2}_m(1) K^{-2,-2}_n{}'(\hat y)  }{ l_3},0\right]
   \\
  =   &\begin{cases}
     [0,0,0,0], & n\ge 1,\, m\in\{1,3,4,5,6,\cdots \} ,\\
    [0 ,   0,     \frac{2K^{-2,-2}_n{}'(\hat y)  }{ l_3 },0], & n\ge 1,\, m=2,
   \end{cases}
 \end{split}
  \\
\label{lem-03-4-curl}
 \nabla \times \bm {q}_{{m,n}}^\ast \big|_{\Gamma}= &\,
  J^{-1}_K(\hat x,\hat y)   K^{-2,-2}_m{}'(\hat x)  K^{-2,-2}_n{}'(\hat y)\big|_{\Gamma}
  =  -\nabla \times \bm {q}_{{m,n}} \big|_{\Gamma}.
 \end{align}
 \end{lemma}


Based on Lemma \ref{lemma-01-4} - Lemma \ref{lemma-03-4},  we shall show in the subsequent  subsections how to construct our $H(\curl^2)$-conforming quadrilateral spectral elements
through  basis functions  \eqref{bases-01} and \eqref{bases-02}  on the reference square.

\subsection{Construction of interior modes}
The interior modes are readily constructed by looking up the trace properties of $\hat {\bs p}_{m,n}$,  $\hat {\bs q}_{m,n}$ and $\hat {\bs q}^*_{m,n}$ in Lemmas \ref{lemma-01-4}-\ref{lemma-03-4}.

 Indeed, \eqref{lem-01-4-curl} and \eqref{lem-01-4-tan} state that both
$  \nabla \times \bm {p}_{{m,n}} \big|_{\Gamma} $ and  $  \bm {p}_{{m,n}}\cdot\bm \tau \big|_{\Gamma}$
vanish if and only if  $ m,n\ge 2$.   While  \eqref{lem-02-4-tan}  and \eqref{lem-02-4-curl}
(resp. \eqref{lem-03-4-tan}  and \eqref{lem-03-4-curl}) imply that
 $  \bm {q}_{{m,n}} \cdot \bm \tau \big|_{\Gamma} $ and $ \nabla \times \bm {q}_{{m,n}} \big|_{\Gamma}$
(resp.  $  \bm {q}^{\ast}_{{m,n}}\cdot \bm \tau \big|_{\Gamma} $  and $ \nabla \times \bm {q}^{\ast}_{{m,n}} \big|_{\Gamma}$)
vanish if and only if  $ m\in\{2,4,5,6,\cdots\}, n\ge 4$  (resp. $ n\in\{2,4,5,6,\cdots\}, m\ge 4$).

By removing the redundant ones, we obtain the interior modes on $K$,
\begin{align*}
&{\bs\phi}_{m,n}(x,y)  := {\bs p}_{m,n}(x,y)  =\bs B_K^{-\tr} \hat  \nabla [ K^{-1,-1}_m(\hat x)  K^{-1,-1}_n(\hat y)] \circ \Phi_K^{-1}, \quad
m, n\ge 2,\\
& {\bs\psi}_{m,n}(x,y) :=  {\bs q}_{m,n}(x,y)  =\bs B_K^{-\tr}[ \hat \nabla K^{-2,-2}_m(\hat x)  K^{-2,-2}_n(\hat y)] \circ \Phi_K^{-1}, \quad m\in\{2,4,5,6,\cdots\}, n\ge 4,\\
& {\bs\psi}_{m,2}(x,y)  :=  {\bs q}^{\ast}_{m,n}(x,y)   =\bs B_K^{-\tr}[\hat K^{-2,-2}_m(\hat x)  \hat\nabla\hat  K^{-2,-2}_n(\hat y)]\circ \Phi_K^{-1}, \quad m\ge 4, n= 2.
\end{align*}

 \subsection{Construction of edge modes} The edge modes should be divided into two types: function edge modes and curl edge modes. We start with the function edge modes. Thanks to Lemma \ref{lemma-01-4}, the functions
 $$ {\bm \phi}_{0,n}=\bs B_K^{-\tr}  \hat \nabla \big[K^{-1,-1}_0(\hat x) K^{-1,-1}_n(\hat y)\big] \circ \Phi_K^{-1},\ n\ge 2,$$  satisfy the criteria for function edge  modes corresponding to $\Gamma_1$.

 Furthermore, according to Lemma \ref{lemma-01-4}, Lemma \ref{lemma-02-4}, and \eqref{Hermite5-01}, we define
\begin{align*}
&\hat {\tilde {\bs  \phi}}_{0,0} = \hat \nabla \big[K^{-1,-1}_0(\hat x) K^{-1,-1}_1(\hat y)\big]
+ \frac12 \hat \nabla \big[K^{-1,-1}_3(\hat x) K^{-1,-1}_1(\hat y)\big]  + \hat \nabla  K^{-2,-2}_2(\hat x) K^{-2,-2}_2(\hat y)
\\
&\qquad   = \left(\frac{3y(y^2-1)(x^2-1)}{16},   \frac{(x+2)(x-1)^2}{8}\right)^{\tr},
\end{align*}
and then find that
\begin{align*}
& {\tilde {\bs  \phi}}_{0,0}\cdot\bm\tau|_\Gamma=\left[ -\frac{1}{ l_1 }, 0  ,0,0\right], &&\nabla\times  {\tilde {\bs  \phi}}_{0,0}|_\Gamma=\left[0,0,0,0\right].
 \end{align*}
 Hence, $ {\tilde {\bs  \phi}}_{0,0}=\bs B_K^{-\tr}\hat{\tilde{\bs \phi}}_{0,0}\circ \Phi_K^{-1}$ or
 $  {\bs  \phi}_{0,0} =\bs B_K^{-\tr} \hat{\bs \phi}_{0,0}  \circ \Phi_K^{-1}:=\bs B_K^{-\tr}\big( \hat{\tilde{\bs \phi}}_{0,0} -\frac18 \hat{\bs \phi}_{3,3}\big) \circ \Phi_K^{-1}=  {\tilde {\bs  \phi}}_{0,0}-  {\bs  \phi}_{3,3}
 $  is another  function edge basis function for $\Gamma_1$.

 The function edge  modes corresponding to $\Gamma_2,\Gamma_3$ and $\Gamma_4$ can be obtained by
making use of  the geometrical and algebraic  symmetry.


Next, let us concentrate on the construction of curl edge  modes. By Lemma \ref{lemma-03-4}, we
 find   that  $
 \bs q^{\ast}_{1,n} =\bs B_K^{-\tr}  [K^{-2,-2}_1(\hat x) \hat \nabla  K^{-2,-2}_n(\hat y)]\circ \Phi_K^{-1},
\   n\in\{2,4,5,6,\cdots\},
$ satisfy
\begin{align*}
&\bs q^{\ast}_{1,n}\cdot \bm \tau|_{\Gamma}=[0,0,0,0],
&& \nabla\times\bs q^{\ast}_{1,n}|_{\Gamma}=\left[\frac{1}{J_K(-1,\hat y)}K_n^{-2,-2}{}'(\hat y),0,0,0\right].
\end{align*}
Since $$J_K|_{\Gamma_1}=\frac{1}{8}\left[l_2l_1\sin\theta_1(1-\hat y)+l_1l_4\sin\theta_4(1+\hat y)\right],$$
we observe that the  curl trace of $\bs q^{\ast}_{1,n}$ on $\Gamma_1$ relies on the geometric quantities $l_1,l_2,l_4,\theta_1,$ and $\theta_4$. While the traces,  up to first order,  of a typical $H(\tc^2)$-conforming basis functions on $\Gamma_i$ only rely on geometric quantities of the edge $\Gamma_i$. This  motivates us to add a multiplier\ $J_K(\hat x,\hat y)$\ to the above basis functions, and define
\begin{align*}
&\bs{\hat  \psi}_{1,n} = J_K(\hat x,\hat y)  K^{-2,-2}_1(\hat x) \hat \nabla  K^{-2,-2}_n(\hat y),
\quad n\in \{2,4,5,6,\cdots\}.
\end{align*}
which finally lead to all curl edge  modes $\bs{ \psi}_{1,n} =\bs B_K^{-\tr}   \bs{\hat  \psi}_{1,n} \circ \Phi_K^{-1}, \ n\in \{2,4,5,6,\cdots\}$, corresponding to $\Gamma_1$ on $K$ such that
\begin{align*}
&\bm \psi_{1,n}\cdot \bm \tau|_{\Gamma}=[0,0,0,0],
&&\nabla\times\bm \psi_{1,n}|_{\Gamma}=[K_n^{-2,-2}{}'(\hat y),0,0,0], \quad n\ \in \{2,4,5,6,\cdots\}.
\end{align*}

 Similarly,  we can obtain other curl edge modes directly by Lemma \ref{lemma-02-4}-Lemma \ref{lemma-03-4}.

\subsection{Construction of vertex modes}
In view of  \eqref{Qfunc1}-\eqref{Qfunc4}, the tangent traces of  all function edge modes on $K$ form a complete system
in $L^2(\partial K)$.   To make  the curl traces of all  $H(\curl^2)$-conforming  basis functions on $K$  a complete system in $C^1(\partial K)$ in regard to \eqref{gamma1curl1}-\eqref{gamma4curl1},
we still need four  vertex modes  whose tangent traces and  curl traces  along $\Gamma$  are zero and  piecewise linear hat functions
along $\partial K$, respectivley.

Suppose ${\bm \psi}_{0,0} = \bs B_K^{-\tr}\hat {\bm \psi}_{0,0}\circ \Phi_K^{-1}$ with
\begin{align}\label{eq-def-01}
\hat {\bm\psi}_{0,0}=  \sum_{1\le m\le 3}  \sum_{n\in \{1,3 \}} \left(a_{m,n} \hat {\bs q}_{m,n} + b_{m,n} \hat {\bs q}^{\ast}_{n,m}\right),
\end{align}
 is such a basis function at  the vertex $P_1$.
It is then expected  that
\begin{align*}
& {\bm\psi}_{0,0}\cdot \bm \tau|_\Gamma=[0,0,0,0],
\qquad
\nabla\times\bm\psi_{0,0}=\left[\frac{1-\hat y}{2},\frac{1-\hat x}{2},0,0\right].
\end{align*}
Indeed,  above requirements can be fulfilled  if $a_{m,n}$ and $b_{m,n}$ satisfy  the following equation,
\begin{align}\label{coeff-04}
\begin{split}
\Big[\frac{1-\hat y}{2},&\frac{1-\hat x}{2},0,0\Big] = \nabla\times\bm \psi_{0,0}|_{\Gamma}
\\
=\bigg[&
-\frac{2( a_{1,1}- b_{1,1})(1-\hat y)(1+3\hat y) + 2(b_{3,1}-a_{1,3}) (1+\hat y)(1-3\hat y) +6b_{2,1}(1-\hat y^2)  }{\left(\left(l_{2} s_{1}-l_{4} s_{4}\right) \hat y-l_{2} s_{1}-l_{4} s_{4}\right) l_{1}},
\\
&\frac{  2(b_{1,1}-a_{1,1}) (1-\hat x) (1+3\hat x)  + 2(b_{1,3}-a_{3,1}) (1+\hat x) (1-3\hat x) + 6a_{2,1}(1-\hat x^2)  }{\left(\left(l_{1} s_{1}-l_{3} s_{2}\right) \hat x-l_{1} s_{1}-l_{3} s_{2}\right) l_{2}},
\\
&-\frac{ 2(a_{3,1}- b_{1,3}) (1-\hat y)(1+3\hat y)    +2(a_{3,3}- b_{3,3}) (1+\hat y)(1-3\hat y)  + 6b_{2,3} (1-\hat y^{2}) }{\left(\left(l_{2} s_{2}-l_{4} s_{3}\right) \hat y-l_{2} s_{2}-l_{4} s_{3}\right) l_{3}},
\\
&\frac{ 2 (b_{3,1}- a_{1,3}) (1-\hat x)(1+3\hat x) +2 (b_{3,3}- a_{3,3})  (1+\hat x)(1-3\hat x)+ 6a_{2,3} (1-\hat x^2) }{\left(\left(l_{1} s_{4}-l_{3} s_{3}\right) \hat x-l_{1} s_{4}-l_{3} s_{3}\right) l_{4}}\bigg],
\end{split}
\end{align}
where  the second equality is established by using  \eqref{lem-02-4-curl} and \eqref{lem-03-4-curl}.

Since the third and fourth entries in \eqref{coeff-04} are zeros,  we  obtain
\begin{align}
\label{coeff_sol}
  b_{3,1}=a_{1,3}, \quad a_{3,1}= b_{1,3}, \quad  b_{3,3}=a_{3,3}, \quad a_{2,3} =b_{2,3}=0,
\end{align}
and  further simplify \eqref{coeff-04}  as
\begin{align}\label{coeff-04a}
\begin{split}
\Big[\frac{1-\hat y}{2},\frac{1-\hat x}{2},0,0\Big]
=\bigg[& \frac{1-\hat y}{2}
\frac{ 12(b_{1,1}-a_{1,1}  -b_{2,1}) \hat y  + 4(b_{1,1}-a_{1,1}-3b_{2,1})     }{\left(\left(l_{2} s_{1}-l_{4} s_{4}\right) \hat y-l_{2} s_{1}-l_{4} s_{4}\right) l_{1}},
\\&
\frac{1-\hat x}{2}\frac{  12(b_{1,1}-a_{1,1}+a_{2,1}) \hat x + 4(b_{1,1}-a_{1,1}+3a_{2,1})  }{\left(\left(l_{1} s_{1}-l_{3} s_{2}\right) \hat x-l_{1} s_{1}-l_{3} s_{2}\right) l_{2}},0,0\bigg].
\end{split}
\end{align}
This implies that
\begin{align}
\label{coeff_sola}
b_{1,1}-a_{1,1}=\frac{l_1l_2s_1}{4}, \quad a_{2,1}=-\frac{1}{6} s_{1} l_{1} l_{2}-\frac{1}{12} s_{2} l_{3} l_{2}, \quad b_{2,1}=\frac{1}{6} s_{1} l_{1} l_{2}+\frac{1}{12} l_{4} l_{1} s_{4}.
\end{align}
Finally,
by selecting
\begin{align}
\label{coeff_solb}
&a_{1,1}=-\frac{s_{1} l_{1} l_{2}}{8},
\quad  b_{1,1}=\frac{s_{1} l_{1} l_{2}}{8},  \quad a_{1,3}=a_{3,1}=a_{3,3}=0,
\end{align}
we get
\begin{align*}
\begin{split}
\hat {\bm\psi}_{0,0}=\,  -\frac{s_{1} l_{1} l_{2}}{8} \hat{\bs q}_{1,1}-\left(\frac{1}{6} s_{1} l_{1} l_{2}+\frac{1}{12} s_{2} l_{3} l_{2}\right)  \hat{\bs q}_{2,1}
 +\frac{s_{1} l_{1} l_{2}}{8} \hat{\bs q}^{\ast}_{1,1} +\left(\frac{1}{6} s_{1} l_{1} l_{2}+\frac{1}{12} l_{4} l_{1} s_{4}\right)
     \hat{\bs q}^{\ast}_{1,2},
   \end{split}
\end{align*}
which yields the vertex mode \eqref{Qvert1}.

By a parity analysis, we can also obtain the  vertex modes \eqref{Qvert2}-\eqref{Qvert4}.

%


\section{Applications to the quad-curl  problem and its eigenvalue problem}
In this section, we propose the $H(\tc^2)$-conforming  quadrilateral spectral element method   to solve the quad-curl problem,
\begin{equation}\label{prob1}
\begin{split}
(\nabla\times)^4\bm u&=\bm f,\quad \text{in}\;\Omega,\\
\nabla \cdot \bm u &= 0, \quad \text{in}\;\Omega,\\
\bm u\times\bm n&=0,\quad \text{on}\;\partial \Omega,\\
\nabla \times \bm u&=0,\quad \text{on}\;\partial \Omega,
\end{split}
\end{equation}
where  $\Omega \in\mathbb{R}^2$ is   Lipschitz domain and  $\bm n$ is the unit outward normal vector to $\partial \Omega$.

In order to satisfy divergence-free condition, we adopt mixed methods  where the constraint $\nabla\cdot\bm u=0$ in \eqref{prob1} is satisfied in a weak sense by introducing an auxiliary unknown $p$.  Hence we adopt the following  variational formulation:
Find $(\bm u; p)\in H_0(\tc^2;\Omega)\times H^1_0(\Omega)$,  s.t.
\begin{align}\label{prob22}
\begin{cases}
\displaystyle a(\bm u,\bm v) + b(\bm v,p)=(\bm f, \bm v), &\forall \bm v\in H_0(\tc^2;\Omega),\\
\displaystyle b(\bm u,q)=0,   &\forall q\in H^1_0(\Omega),
\end{cases}
\end{align}
where
\begin{align*}
&a(\bm u,\bm v)=((\nabla\times)^2 \bm u, (\nabla\times)^2 \bm v), \\
&b(\bm v,p)=(\bm v,\nabla p ).
\end{align*}
The well-posedness of the variational problem can be found in \cite{Sun2019Awg}.

\subsection{Approximation spaces}

 Let $L,M,N$ be  three  integers $\ge 3$. 
We  introduce the following mapped polynomial space,
\begin{align}
\label{VLMN}
\begin{split}
V_{L,M,N}(K)=\Big\{&\bm \phi_{m,n}, 2\le m,n\le L, \bm \phi_{m,0}, \bm \phi_{1,n}, 1\le m,n \le M,\\
&\bm \phi_{m,1}, \bm \phi_{0,n},  m,n\in\{0,2,3,4,\cdots,M\},\\
&\bm \psi_{i,j}, 0\le i,j\le 1,\\
&\bm \psi_{m,n}, \bm \psi_{m,2},\bm \psi_{2,n},  4\le m, n\le N,\\
&\bm \psi_{1,l}, \bm \psi_{l,1},\bm \psi_{3,l}, \bm \psi_{l,3}, l\in\{2,4,5,6,\cdots,N\}\Big\}.
\end{split}
\end{align}
We shall abbreviate   $V_{N,N,N}(K)$ as  $V_{N}(K)$.

We now list all the $H(\tc^2)$-conforming basis functions in $V_{L,M,N}(K)$  in Table \ref{basis-ref1},
which shows that $\displaystyle \min_{L,M, N \ge 3}  \mathrm{dim} V_{L,M, N}= 24$.
\begin{table}[!htbp]
\centering
\caption{The $H(\tc^2)$-conforming elements  on the element $ K$.}\label{basis-ref1}
\begin{tabular}{lllll}
  \hline
  modes & basis & cardinality \\
  \hline
  \multirow{3}{0cm}{interior} & ${\bm \phi}_{m,n}, 2\le m,n\le L$   &\multirow{3}{2.5cm}{ \small{$  (L-1)^2+(N-1)(N-3)$}}\\
  &${\bm \psi}_{m,n}, m\in \{2,4,5,6,\cdots,N\},4 \le n \le N$ &\\
  &${\bm \psi}_{m,2},  4\le m\le N$& \\
  \hline
  \multirow{4}{2cm}{function edge}
  &$\Gamma_4$: ${\bm \phi}_{m,1}, m\in\{0,2,3,4\cdots,M\}$&{\multirow{4}{1cm}{$4M$}}\\
  &$\Gamma_3$: ${\bm \phi}_{1,n}, 1\le n \le M$&\\
  &$\Gamma_2$: ${\bm \phi}_{m,0}, 1\le m \le M$&\\
  &$\Gamma_1$: ${\bm \phi}_{0,n}, n\in\{0,2,3,4\cdots,M\}$& \\
  \hline
  \multirow{4}{2cm}{curl edge} & $\Gamma_4: {\bm \psi}_{m,3}, m\in\{2,4,5,6,\cdots,N\}$&\multirow{3}{1.5cm}{$4(N-2)$}\\
  &$\Gamma_3: {\bm \psi}_{3,n}, n\in\{2,4,5,6,\cdots,N\}$&\\
  &$\Gamma_2: {\bm \psi}_{m,1}, m\in\{2,4,5,6,\cdots,N\}$  & \\
  &$\Gamma_1: {\bm \psi}_{1,n},  n\in\{2,4,5,6,\cdots,N\}$  & \\
  \hline
  {vertex} &${\bm\psi}_{i,j}, 0\le i,j\le 1$& ${4}$ \\
  \hline
\end{tabular}
\end{table}

 The following lemma states the polynomial spaces $V_{L, M,N}(K),\ L, M, N \ge 3$ contains lower-order polynomials on arbitrary quadrilateral $K$, thus  form a complete system  in $H(\curl^2; K)$.

\begin{lemma}\label{lowerorderpolynomials}
It holds that
\begin{align}
(1,0)^{\tr}=& -x_{14} {\bm \phi}_{0,0}+x_{21} {\bm \phi}_{1,0}+x_{32} {\bm \phi}_{1,1}-x_{43} {\bm \phi}_{0,1},\label{low-01}
\\
(0,1)^{\tr}=& -y_{14} {\bm \phi}_{0,0}+y_{21} {\bm \phi}_{1,0}+y_{32} {\bm \phi}_{1,1}-y_{43} {\bm \phi}_{0,1},\label{low-02}
\\
\begin{split}
(x,0)^{\tr}=&  -\frac{(x_1^2-x_4^2)}{2}  {\bm \phi}_{0,0}+\frac{(x_2^2-x_1^2)}{2} {\bm \phi}_{1,0}+\frac{(x_3^2-x_2^2)}{2} {\bm \phi}_{1,1}-\frac{(x_4^2-x_3^2)}{2}  {\bm \phi}_{0,1}
\\
&+ \frac{x_{14}^2}{2} \bm \phi_{0,2}+\frac{x_{21}^2}{2} \bm \phi_{2,0}+\frac{x_{32}^2}{2}  \bm \phi_{1,2}+\frac{x_{43}^2}{2} \bm \phi_{2,1}+\frac{(x_{32}+x_{14})^2}{2}  \bm \phi_{2,2} ,
\end{split}\label{low-03}
\\
\begin{split}
(0,x)^{\tr}=&  -\frac{(x_1+x_4)y_{14}}{2}{\bm \phi}_{0,0}+\frac{(x_2+x_1)y_{21}}{2}  {\bm \phi}_{1,0}+\frac{(x_3+x_2)y_{32}}{2} {\bm \phi}_{1,1}
\\
&-\frac{(x_4+x_3)y_{43}}{2} {\bm \phi}_{0,1}+ \frac{x_{14}y_{14}}{2}  \bm \phi_{0,2}+\frac{x_{21}y_{21}}{2}  \bm \phi_{2,0}+\frac{x_{32}y_{32}}{2} \bm \phi_{1,2}
\\
&+\frac{x_{43}y_{43}}{2} \bm \phi_{2,1}
+\frac{(x_{32}+x_{14})(y_{32}+y_{14})}{2}  \bm \phi_{2,2}+(  \tilde{\bm \psi}_{0,0}+ \tilde{\bm \psi}_{0,1}+ \tilde{\bm \psi}_{1,0}+ \tilde{\bm \psi}_{1,1}) ,
\end{split}
\label{low-04}
\\
\begin{split}
(0,y)^{\tr}=&  -\frac{(y_1^2-y_4^2)}{2}  {\bm \phi}_{0,0}+\frac{(y_2^2-y_1^2)}{2} {\bm \phi}_{1,0}+\frac{(y_3^2-y_2^2)}{2}  {\bm \phi}_{1,1}-\frac{(y_4^2-y_3^2)}{2}{\bm \phi}_{0,1}
\\
&+ \frac{y_{14}^2}{2}  \bm \phi_{0,2}+\frac{y_{21}^2}{2}  \bm \phi_{2,0}+\frac{y_{32}^2}{2} \bm \phi_{1,2}+\frac{y_{43}^2}{2}  \bm \phi_{2,1}+\frac{(y_{32}+y_{14})^2}{2}  \bm \phi_{2,2} ,
\end{split}
\label{low-05}
\\
\begin{split}
(y,0)^{\tr}=&  -\frac{x_{14}(y_1+y_4)}{2}{\bm \phi}_{0,0}+\frac{x_{21}(y_2+y_1)}{2} {\bm \phi}_{1,0}+\frac{x_{32}(y_3+y_2)}{2} {\bm \phi}_{1,1}
\\
&-\frac{x_{43}(y_4+y_3)}{2}  {\bm \phi}_{0,1}+ \frac{x_{14}y_{14}}{2} \bm \phi_{0,2}+\frac{x_{21}y_{21}}{2} \bm \phi_{2,0}+\frac{x_{32}y_{32}}{2} \bm \phi_{1,2}
\\
&+\frac{x_{43}y_{43}}{2} \bm \phi_{2,1}
+\frac{(x_{32}+x_{14})(y_{32}+y_{14})}{2} \bm \phi_{2,2}-( \tilde{\bm \psi}_{0,0}+ \tilde{\bm \psi}_{0,1}+ \tilde{\bm \psi}_{1,0}+ \tilde {\bm \psi}_{1,1}),
\end{split}
\label{low-06}\\
\begin{split}
1=&\nabla\times\Big(  -\frac{(x_1+x_4)y_{14}}{2}  {\bm \phi}_{0,0}+\frac{(x_2+x_1)y_{21}}{2} {\bm \phi}_{1,0}+\frac{(x_3+x_2)y_{32}}{2}  {\bm \phi}_{1,1}
\\
&-\frac{(x_4+x_3)y_{43}}{2}  {\bm \phi}_{0,1}+(\tilde{\bm \psi}_{0,0}+ \tilde{\bm \psi}_{0,1}+ \tilde{\bm \psi}_{1,0}+\tilde {\bm \psi}_{1,1}) \Big),
\end{split}\label{low-07}
\end{align}
where
\begin{align*}
 {\tilde{\bm\psi}}_{0,0}= {{\bm\psi}}_{0,0}&-\frac{1}{48}\left(2s_1l_1l_2+s_4l_4 l_1\right) {\bm \phi}_{2,3}+\frac{1}{48}\left(2s_1l_1l_2+s_2l_2l_3\right) {\bm \phi}_{3,2},  \\
 {\tilde{\bm\psi}}_{0,1}= {{\bm\psi}}_{0,1}&+\frac{1}{48}\left(2s_2l_2l_3+s_1l_1l_2\right) {\bm \phi}_{3,2}+\frac{1}{48}\left(2s_2l_2l_3+s_3l_3l_4\right) {\bm \phi}_{2,3},  \\
 {\tilde{\bm\psi}}_{1,1}= {{\bm\psi}}_{1,1}&+\frac{1}{48}\left(2s_3l_3l_4+s_2l_2l_3\right) {\bm \phi}_{2,3}-\frac{1}{48}\left(2s_3l_3l_4+s_4l_4l_1\right) {\bm \phi}_{3,2},\\
 {\tilde{\bm\psi}}_{1,0}= {{\bm\psi}}_{1,0}&-\frac{1}{48}\left(2s_4l_4l_1+s_3l_3l_4\right) {\bm \phi}_{3,2}-\frac{1}{48}\left(2s_4l_4l_1+s_1l_1l_2\right) {\bm \phi}_{2,3},
\end{align*}
are alternative vertex modes.
\end{lemma}
The proof is postponed to Appendix \ref{AppendixA}.

 Based on Lemma \ref{lowerorderpolynomials}, we can    further define two simple approximation spaces $V_1(K)=V_{1,1,1}(K)$ and $V_2(K)=V_{2,2,2}(K)$ on $K$  with
the lowest  DOFs:
\begin{align}
\label{V1}
&V_{1,1,1}(K)=:\text{span}\big\{ {\bm\phi}_{0,0},  {\bm\phi}_{1,0}, {\bm \phi}_{1,1},  {\bm\phi}_{0,1},    \tilde{\bm \psi}_{0,0}, \tilde{ \bm\psi}_{0,1}, \tilde{\bm\psi}_{1,1}, \tilde{\bm\psi}_{1,0}\big\},
\\
\label{V2}
 &V_{2,2,2}(K):=\text{span}\big\{{\bm\phi}_{0,0},  {\bm\phi}_{1,0},  {\bm\phi}_{1,1},  {\bm\phi}_{0,1}
,  \bm\phi_{0,2}, \bm\phi_{2,0},  \bm\phi_{1,2},  \bm\phi_{2,1},  \bm\phi_{2,2}, \tilde{\bm\psi}_{0,0},  \tilde{\bm\psi}_{0,1},  \tilde{\bm\psi}_{1,1},  \tilde{\bm\psi}_{1,0}\big\},
\end{align}
such that  $(1,0), (0,1) \in V_1(K)$;   $(1,0), (0,1), (x,0),  (0,x), (0,y), (y,0) \in V_2(K)$; and
$1\in \curl  V_1(K) \subset \curl   V_2(K)$.
It is obvious that $\dim V_1(K)=8$ and $ \dim V_2(K)\le 13$.  Whenever  $K$  is a  parallelogram, $\dim V_2(K)=12$.

\subsection{Approximation scheme}
Let \,$\mathcal{T}_h=\{K_i\}\,$ be a partition of the domain $\Omega$ of the mesh size $h$
consisting of convex quadrilaterals.
We assume that $\mathcal{T}_h$ is regular, i.e.,  the intersection
$K_i\cap  K_j,\ i \neq j$ is either empty or a node or an entire edge of both $\bar{K}_i$ and $\bar{K}_j$.

Let $L,M,N$ be an integer triplet such that  $L,M,N\ge 3$ or $1\le L=M=N\le 2$ hereafter.  We define the $H(\curl^2)$-conforming approximation spaces for  the vector filed $\bs u$,
\begin{eqnarray*}
	&&  W_{L,M,N}^h=\left\{\bm{v}_{L,M,N}^h\in H(\text{curl}^2;\Omega):\ \bm{v}_{L,M,N}^h|_K\in V_{L,M,N}(K),\ \forall K\in\mathcal{T}_h\right\},\\
	&&   \OO{W}_{L,M,N}^h=\left\{\bm{v}_{L,M,N}^h\in W_{L,M,N}^h:\ \bm{n} \times \bm{v}_{L,M,N}^h=0\ \text{and}\ \nabla\times  \bm{v}_{L,M,N}^h = 0 \ \text {on} \ \partial\Omega\right\}.
\end{eqnarray*}
Moreover, we introduce the scalar functions
\begin{align*}
\varphi_{m,n} = \hat \varphi_{m,n}\circ  \Phi_K^{-1}   \text{ on }  K, \qquad \hat \varphi_{m,n}(x,y)=[K_m^{-1,-1}( \hat x)K_n^{-1,-1}( \hat y)]  \text{ on } \hat K,
\end{align*}
and the corresponding local function space
\begin{align*}
R_{L,M}(K):=\big\{ \hat \varphi_{m,n}\circ  \Phi_K^{-1} : 2\le m,n\le L\big\} \oplus  \big\{\hat \varphi_{m,n}\circ  \Phi_K^{-1}:   0 \le   m,n  \le M; \, \min(m,n) \le 1 \big \},
\end{align*}
such that $ \left\{ \bs u\in V_{L,M,N}(K):  \nabla \times \bs u =0\right\} = \left\{ \nabla u: u\in R_{L,M}(K) \right\}$.
Indeed,
by \eqref{Qfunc1}-\eqref{Qfunc4},
\begin{align*}
&\hat {\bs\phi}_{1,1} + \hat {\bs\phi}_{0,1} = \hat \nabla\hat \varphi_{1,1}, && \hat {\bs\phi}_{0,0} - \hat {\bs\phi}_{0,1} =\hat  \nabla \hat \varphi_{0,1}, \\
& -\hat {\bs\phi}_{0,0} - \hat {\bs\phi}_{1,0} = \hat \nabla \hat \varphi_{0,0}, && -\hat {\bs\phi}_{1,1} + \hat {\bs\phi}_{1,0} =\hat  \nabla \hat \varphi_{1,0},
\end{align*}
which together with \eqref{Qinterior}  and \eqref{Qfunc1}-\eqref{Qfunc4} implies
\begin{align*}
  \nabla \varphi=  [\bs B_K^{-\tr}\hat \nabla {\hat \varphi}]\circ  \Phi_K^{-1}  \in  V_{L,M,N}(K), \qquad \forall \varphi=\hat \varphi\circ  \Phi_K^{-1}\in R_{L,M}(K).
\end{align*}
We now define the  $H^1$-conforming approximation spaces for the auxiliary  function $p$,

\begin{eqnarray*}
	&&  S_{L,M}^h=\left\{{w}_{L,M}^h\in H^1(\Omega):\  w_{L,M}^h|_K\in R_{L,M}(K)\right\},\\
	&&  \OO{S}_{L, M}^h=\left\{{w}_{L,M}^h\in S_{L,M}^h,\;{w}_{L,M}^h=0\ \text {on} \  {\partial\Omega}\right\}.
\end{eqnarray*}
Once again, we shall abbreviate   $W_{N,N,N}^h$  and $\OO{W}_{N,N,N}^h$ as  $W_{N}^h$ and $\OO{W}_{N}^h$, respectively;
and abbreviate $S^h_{N,N} $ and $\OO{S}^h_{N,N} $ as  $S_{N}^h$ and $\OO{S}_{N}^h$, respectively.


The $H(\tc^2)$-conforming quadrilateral spectral element method for \eqref{prob22} seeks $(\bm u^h_{L,M,N};p^h_{L,M})\in \OO{W}^{h}_{L,M,N}\times \OO{S}^{h}_{L,M}$,  s.t.
\begin{equation}\label{prob3}
\begin{cases}
\displaystyle a(\bm u^h_{L,M,N},\bm v^h_{L,M,N}) + b(\bm v^h_{L,M,N},p_{L,M}^h)=(\bm f, \bm v^h_{L,M,N}), & \displaystyle  \forall \bm v^h_{L,M,N}\in \OO{W}^{h}_{L,M,N},\\
\displaystyle b(\bm u^h_{L,M,N},q_{L,M}^h)=0, & \displaystyle  \forall q^h_{L,M}\in \OO{S}^{h}_{L,M}.
\end{cases}
\end{equation}


It  is obvious that
\begin{align}
\label{SW}
\nabla \OO{S}_{ L,M}^h \subset \OO{W}_{L,M,N}^h.
\end{align}
As a result,
\begin{align*}
&\sup_{ \bs u \in \OO{W}_{L,M,N}^h } \frac{(\bs u, \nabla p) }{\|\bs u\|_{H(\curl^2;\Omega)}}
\overset{\bs u = \nabla p}\ge  \frac{(\nabla p, \nabla p) }{\|\nabla p\|_{H(\curl^2;\Omega)}}
=  \frac{\|\nabla p\|_{L^2(\Omega)}^2 }{\|\nabla p\|_{L^2(\Omega)}} = \|\nabla p\|_{L^2(\Omega)}, \quad p\in \OO{S}_{L, M}^h.
\end{align*}
Meanwhile, it is easy to see that
\begin{align*}
&\frac{(\bs u, \nabla p) }{\|\bs u\|_{H(\curl^2;\Omega)}}
\le   \frac{\|\bs u\|_{L^2(\Omega)} \, \| \nabla p\|_{L^2(\Omega)}  }{\|\bs u\|_{H(\curl^2;\Omega)}}
\le \| \nabla p\|_{L^2(\Omega)}.
\end{align*}
Consequently,
\begin{align*}
\sup_{ \bs u \in \OO{W}_{L,M,N}^h }\frac{(\bs u, \nabla p) }{\|\bs u\|_{H(\curl^2;\Omega)}} = \|\nabla p\|_{L^2(\Omega)},
\end{align*}
which shows that the discrete Ladyzhenskaya-Babu\v{s}ka-Brezzi condition is satisfied, thus \eqref{prob3}
is well-posed.

Assembling the global ``stiffness'' matrix and  ``damping'' matrix $\bs A, \bs B$, we arrive at the following equivalent
algebraic system,
\begin{align}\label{alg-source}
\begin{pmatrix}
\bs A &\bs B\\
\bs B^{\tr} & \bs 0\\
\end{pmatrix}\begin{pmatrix}
  \mathfrak{u}\\
  \mathfrak{p}\\
\end{pmatrix}=\begin{pmatrix}
\bs F\\
0\\
\end{pmatrix},
\end{align}
where $ \mathfrak{ u},   \mathfrak{p}$ are two column vectors which represent  the DOFs corresponding to the global basis functions, respectively.

\section{Numerical result}
\subsection{The source problem}
In this subsection, we consider the problem \eqref{prob1} on a unit square $\Omega=(0,1)^2$ with exact solution
	\begin{equation}\label{exactsolution}
	\bm u=\left(
	\begin{array}{c}
	3\pi\sin^3(\pi x)\sin^2(\pi y)\cos(\pi y) \\
	-3\pi \sin^3(\pi y)\sin^2(\pi x)\cos(\pi x)\\
	\end{array}
	\right).
	\end{equation}
	The source term $\bm f$ can be obtained by a simple calculation. 
	We denote
	\[\bm e_{L,M,N}^h=\bm u-\bm u_{L,M,N}^h.\]
and simply write $\bm e_{N,N,N}^h$ as $\bm e_{N}^h$.

\begin{figure}[htbp]
  \subfigure[$n=1$;]{\includegraphics[width=0.32\textwidth]{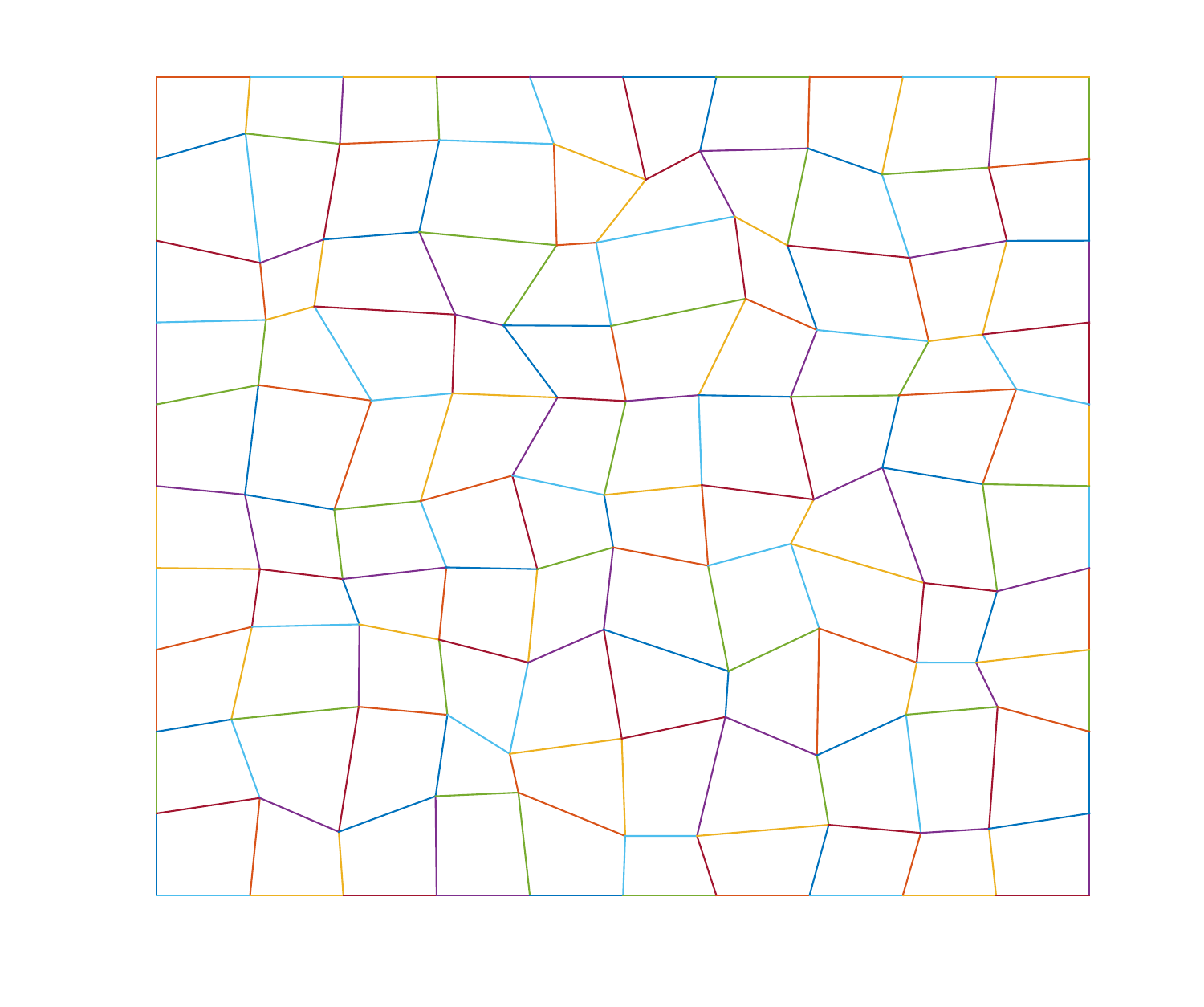}}
  \subfigure[$n=2$;]{\includegraphics[width=0.32\textwidth]{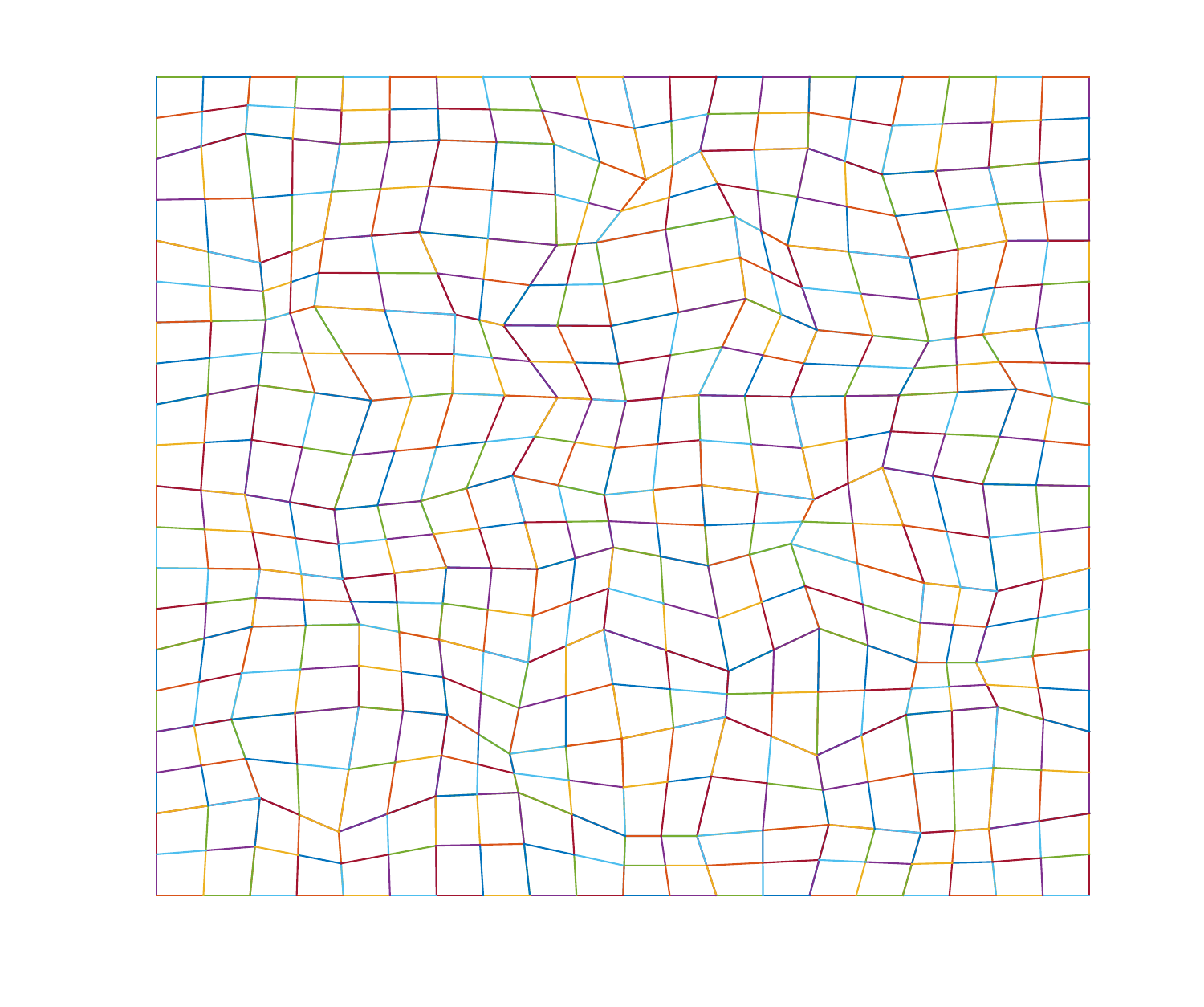}}
  \subfigure[$n=3$.]{\includegraphics[width=0.32\textwidth]{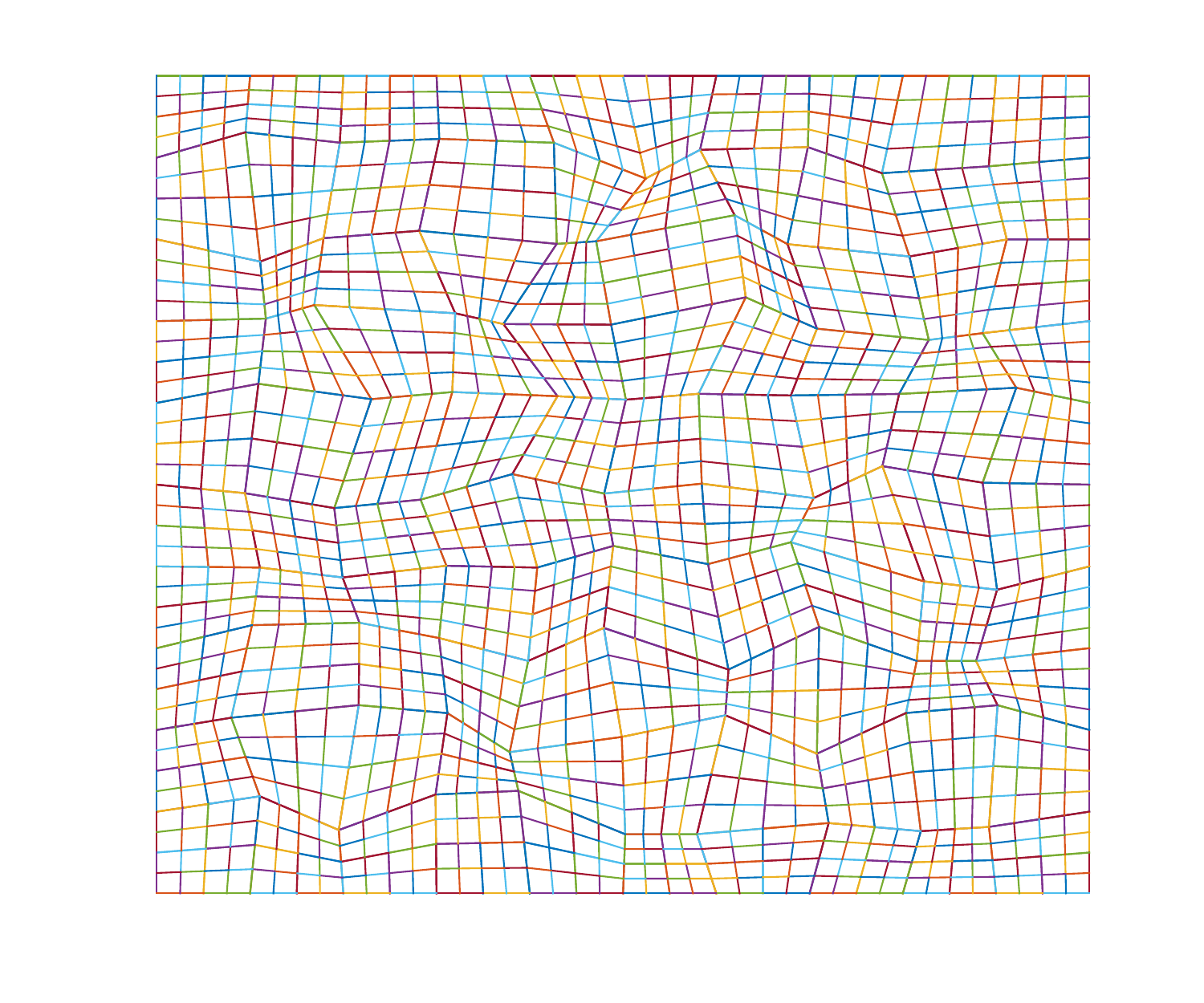}}
  \caption{Non-uniform quadrilateral meshes.}\label{quadrilateral_mesh}
\end{figure}

Various tests are implemented to demonstrate the  validity  and efficiency of the $H(\tc^2)$-conforming  spectral element method.
 We begin with the $h$-convergence  of the  simplest  two approximation spaces  $\OO{W}^h_1$ and $ \OO{W}^h_2$.  The initial partition is a nonuniform quadrilateral mesh  with $h=10^{-1}$ (see Figure \ref{quadrilateral_mesh} (a)),  followed by several levels of subsequent meshes using the regular refinement (see Figure \ref{quadrilateral_mesh} (b) and (c) for the first two  levels).
Approximation errors in $L^2$-, $H(\tc)$- and $H(\tc^2)$-seminorms, i.e., $\|\bm e_{N}^h\|$, $\|\nabla \times\bm e_{N}^h\|$ and $\|\nabla \times \nabla \times \bm e_{N}^h\|$,
are obtained in  Table \ref{tablower1} and  Table \ref{tablower2} for $N=1$ and  $N =2$, respectively.

With the approximation space $ \OO{W}_{1}^h$, the simplest  $H(\tc^2)$-conforming quadrilateral  method   consists 8  DOFs on each physical element.
The vector field  $\bs u_1^h$ converges to $\bs u$ at  orders of $\mathcal{O}(h)$,  $\mathcal{O}(h^2)$ and $\mathcal{O}(h)$  in  semi-norms in $L^2(\Omega)$, $H(\tc;\Omega)$ and $H(\tc^2;\Omega)$, respectively.  This situation agrees well with the convergence  behavior   of the simplest rectangular elements reported  in  \cite{Hu2020simple}.

\begin{table}[h!]
		\centering
		\caption{Numerical results by using  $H(\tc^2)$-conforming  elements with $L=M=N=1$.} \label{tablower1}
		\begin{tabular}{cccccccc}
			\hline
			$h$ &$\left\|\bm e_{1}^h\right\|$& rates&$\left\|\nabla\times\bm e_{1}^h\right\|$& rates&$\left\|\nabla\times\nabla\times\bm e_{1}^h\right\|$& rates\\
			\hline
			$1/10$&2.7016060e-01&      &8.9911798e-01&      &31.9663395& \\
			$1/20$&1.3087910e-01&1.0456&1.9697265e-01&2.1905&14.7861071&1.1123\\
			$1/40$&6.5062865e-02&1.0083&4.7352656e-02&2.0565&7.23589118&1.0310\\
			$1/80$&3.2494568e-02&1.0016&1.1725121e-02&2.0138&3.59895334&1.0076\\
			$1/160$&1.6243391e-02&1.0003&2.9243540e-03&2.0034&1.79714478&1.0019\\
			  \hline
		\end{tabular}
	\end{table}

 Note that the convergence rate of $\bs u_1^h$ is lower than that of $\nabla \times \bs u_1^h$.  This is very different from convergence behaviors  of elliptic equations with the gradient  operator.

In order to improves the rate of convergence for $\bs u_1^h$, more DOFs are needed. So we carry out numerical experiments
on the $H(\tc^2)$-conforming quadrilateral  method with   the approximation space $ \OO{W}_{2}^h$.     We see from Table \ref{tablower2} that with 13 element DOFs ($L=M=N=2$), the convergence rate  of $\bs u_2^h$ raises to $O(h^2)$, while the convergence rates (even accuracy) of $\nabla \times\bs u_2^h$ and $\nabla \times \nabla \times\bs u_2^h$ are unchanged.

\begin{table}[h!]
		\centering
		\caption{Numerical results by using  $H(\tc^2)$-conforming quadrilateral elements with $L=M=N=2$.} \label{tablower2}
		\begin{tabular}{cccccccc}
			\hline
			$h$ &$\left\|\bm e_{2}^h\right\|$& rates&$\left\|\nabla\times\bm e_{2}^h\right\|$& rates&$\left\|\nabla\times\nabla\times\bm e_{2}^h\right\|$& rates\\
			\hline
			$1/10$&9.1992742e-02&      &8.7767063e-01&      &31.8557504& \\
			$1/20$&2.0157951e-02&2.1902&1.9253484e-01&2.1886&14.7330986&1.1125\\
			$1/40$&4.8913038e-03&2.0431&4.6475991e-02&2.0506&7.21785897&1.0294\\
			$1/80$&1.2142849e-03&2.0101&1.1519875e-02&2.0124&3.59091342&1.0072\\
			$1/160$&3.0305706e-04&2.0024&2.8738828e-03&2.0031&1.79323838&1.0018\\
			  \hline
		\end{tabular}
	\end{table}

Whenever $K$   becomes a rectangle,  the DOFs  reduces to 12, which is one DOF less than
the rectangular  element in  \cite{Hu2020simple}.

Next,  let us examine the $h$-convergence of $\bs u_{L,M,N}^h$  with $L=N=3$ and $M=N,N+1$.
Table \ref{tablower3} and Table \ref{tablower4} show that the errors $\|\bm e_{N,M,N}^h\|$, $\|\nabla \times\bm e_{N,M,N}^h\|$ and $\|\nabla \times \nabla \times \bm e_{N,M,N}^h\|$ decay asymptotically as  $\mathcal{O}(h^M)$,  $\mathcal{O}(h^N)$ and $\mathcal{O}(h^{N-1})$, respectively. It implies that the convergence order  in  the $L^2$-norm will increase by one if
one supplements 4 DOFs on each element to change the approximation space $ \OO{W}^h_{N,N,N}$ to $ \OO{W}^h_{N,N+1,N}$.
Similar convergence behavior of $\bs u_{N,M,N}^h$  can be observed in Table \ref{tab1} and Table \ref{tab2} for
$L=N=4$ and $M=N,N+1$.

\begin{table}[h!]
		\centering
		\caption{Numerical results by  using the  $H(\tc^2)$-conforming elements with $L=M=N=3$.} \label{tablower3}
		\begin{tabular}{cccccccc}
			\hline
			$h$ &$\left\|\bm e_{3}^h\right\|$& rates&$\left\|\nabla\times\bm e_{3}^h\right\|$& rates&$\left\|\nabla\times\nabla\times\bm e_{3}^h\right\|$& rates\\
			\hline
			$1/10$&1.4531621e-02&       &2.356674e-01&      &1.3413460e+01& \\
			$1/20$&1.1122290e-03&3.7076&2.875953e-02&3.0346&3.6870789e+00&1.8631\\
			$1/40$&8.3518345e-05&3.7352&3.672498e-03&2.9692&9.7564651e-01&1.9180\\
			$1/80$&6.6964145e-06&3.6406&4.640792e-04&2.9843&2.5019086e-01&1.9633\\
			$1/160$&6.4506166e-07&3.3759&5.825095e-05&2.9940&6.3218033e-02&1.9846\\
			  \hline
		\end{tabular}
	\end{table}

\begin{table}[h!]
		\centering
		\caption{Numerical results by using the   $H(\tc^2)$-conforming elements with $M=4,\ L=N=3$.} \label{tablower4}
		\begin{tabular}{cccccccc}
			\hline
			$h$ &$\left\|\bm e_{3,4,3}^h\right\|$& rates&$\left\|\nabla\times\bm e_{3,4,3}^h\right\|$& rates&$\left\|\nabla\times\nabla\times\bm e_{3,4,3}^h\right\|$& rates\\
			\hline
			$1/10$&1.4367808e-02&       &2.3566738e-01&      &1.3413460e+01& \\
			$1/20$&1.0749797e-03&3.7405&2.8759528e-02&3.0346&3.6870789e+00&1.8631\\
			$1/40$&7.5421333e-05&3.8332&3.6724978e-03&2.9692&9.7564651e-01&1.9180\\
			$1/80$&4.9693459e-06&3.9238&4.6407944e-04&2.9843&2.5019086e-01&1.9633\\
			  \hline
		\end{tabular}
	\end{table}

\begin{table}[h!]
		\centering
		\caption{Numerical results by using the  $H(\tc^2)$-conforming elements with $L=M=N=4$.} \label{tab1}
		\begin{tabular}{cccccccc}
			\hline
			$h$ &$\left\|\bm e_{4}^h\right\|$& rates&$\left\|\nabla\times\bm e_{4}^h\right\|$& rates&$\left\|\nabla\times\nabla\times\bm e_{4}^h\right\|$& rates\\
			\hline
			$1/10$&3.5791271e-04&      &1.4638912e-02&      &1.3023272& \\
			$1/20$&1.6019024e-05&4.4817&9.5069288e-04&3.9447&1.6959892e-01&2.9409\\
			$1/40$&8.6362784e-07&4.2132&5.9743054e-05&3.9921&2.1604928e-02&2.9727\\
			$1/80$&5.1620409e-08&4.0644&3.7288683e-06&4.0020&2.7236022e-03&2.9878\\
			  \hline
		\end{tabular}
	\end{table}

\begin{table}[h!]
  \centering
  \caption{Numerical results by using the   $H(\tc^2)$-conforming elements with $M=5, L=N=4$.} \label{tab2}
  \begin{tabular}{cccccccc}
   \hline
   $h$ &$\left\|\bm e_{4,5,4}^h\right\|$& rates&$\left\|\nabla\times\bm e_{4,5,4}^h\right\|$& rates&$\left\|\nabla\times\nabla\times\bm e_{4,5,4}^h\right\|$& rates\\
   \hline
   $1/10$&1.8193533e-04&      &1.0483422e-02&      &9.1705319e-01& \\
   $1/20$&5.7613686e-06&4.9809&6.6294347e-04&3.9831&1.1567122e-01&2.9870\\
   $1/40$&1.8301771e-07&4.9764&4.1611226e-05&3.9938&1.4503254e-02&2.9956\\
   $1/80$&5.7630487e-09&4.9890&2.6033668e-06&3.9985&1.8134660e-03&2.9996\\
     \hline
  \end{tabular}
 \end{table}
 We point out that our $H(\curl^2)$-conforming elements  outperform the ones designed in \cite{Hu2020simple},
since our method using approximation spaces $\OO{W}_{N,N+1,N}^h$  ($N\ge 3$) with  $2N^2+2N+4$ DOFs   on each quadrilateral element  provides convergence orders  $\mathcal{O}(h^{N+1})$ in $L^2(\Omega)$,  $\mathcal{O}(h^N)$ in $H(\tc;\Omega)$  and $\mathcal{O}(h^{N-1})$ in $H(\tc^2;\Omega)$
for the numerical vector field, while \cite{Hu2020simple} needs $2N^2+4N+3$ DOFs on each  rectangular element  to
acquire the same orders of convergence.

 Interestingly, the situation  where the convergence rate of $\bs u_{L,M,N}^h$ is lower than that of $\nabla \times \bs u_{L,M,N}^h$
can be reproduced  for any $L=N-1, M=N-1$ with $N\ge 4$. Indeed, Table \ref{tab3-01} shows that the errors $\|\bm e_{N-1,N-1,N}^h\|$, $\|\nabla \times\bm e_{N-1,N-1,N}^h\|$ and $\|\nabla \times \nabla \times \bm e_{N-1,N-1,N}^h\|$ decay asymptotically as  $\mathcal{O}(h^{N-1})$,  $\mathcal{O}(h^N)$ and $\mathcal{O}(h^{N-1})$, respectively.

 \begin{table}[h!]
  \centering
\caption{Numerical results by using the   $H(\tc^2)$-conforming elements with $L=M=3, N=4$.} \label{tab3-01}
  \begin{tabular}{cccccccc}
   \hline
   $h$ &$\left\|\bm e_{3,3,4}^h\right\|$& rates&$\left\|\nabla\times\bm e_{3,3,4}^h\right\|$& rates&$\left\|\nabla\times\nabla\times\bm e_{3,3,4}^h\right\|$& rates\\
   \hline
   $1/10$&1.7929703e-03&      &5.9958211e-03&      &5.6265349e-01& \\
   $1/20$&2.2577801e-04&2.9894&3.8206331e-04&3.9721&7.1249756e-02&2.9813\\
   $1/40$&2.8439099e-05&2.9890&2.4014441e-05&3.9918&8.9231314e-03&2.9973\\
   $1/80$&3.5631885e-06&2.9966&1.5037520e-06&3.9973&1.1151661e-03&3.0003\\
     \hline
  \end{tabular}
 \end{table}

\begin{figure}
  \centering
   \subfigure[$h=1$ for single element;]{ \includegraphics[width=0.32\textwidth]{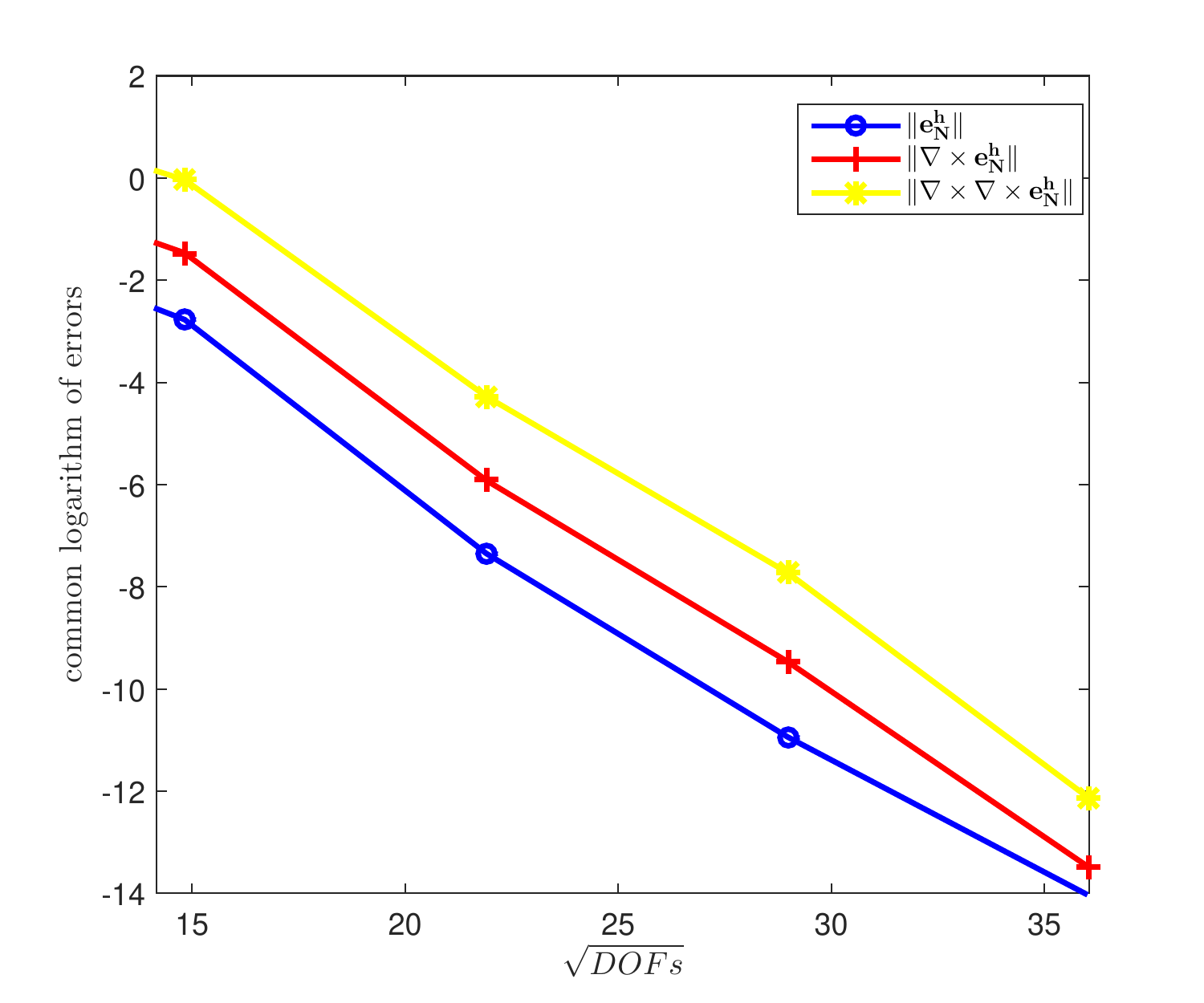}}%
   \subfigure[$h=1/2$ for 4 elements;]{ \includegraphics[width=0.32\textwidth]{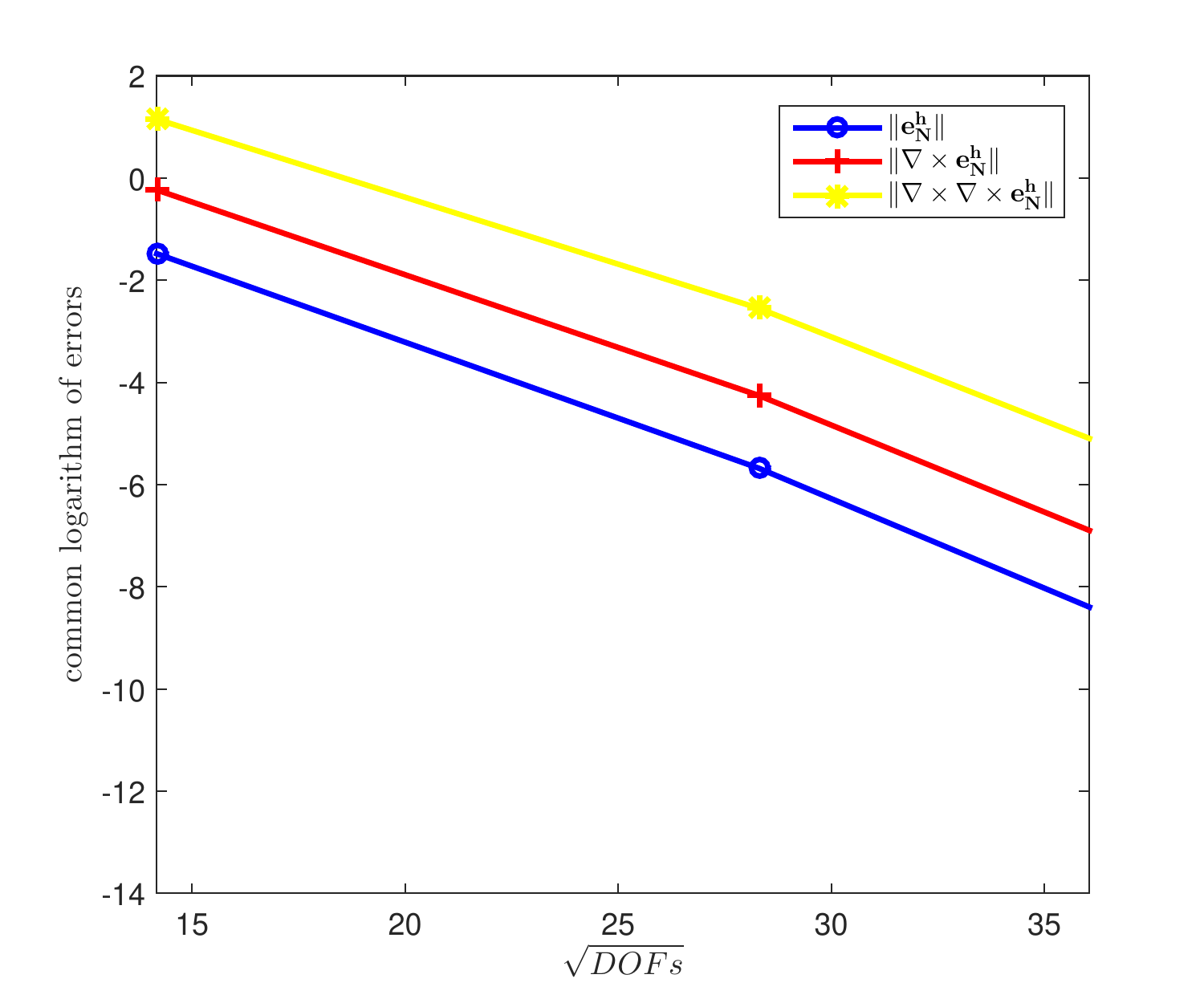}}%
   \subfigure[$h=1/8$ for 64 elements.]{ \includegraphics[width=0.32\textwidth]{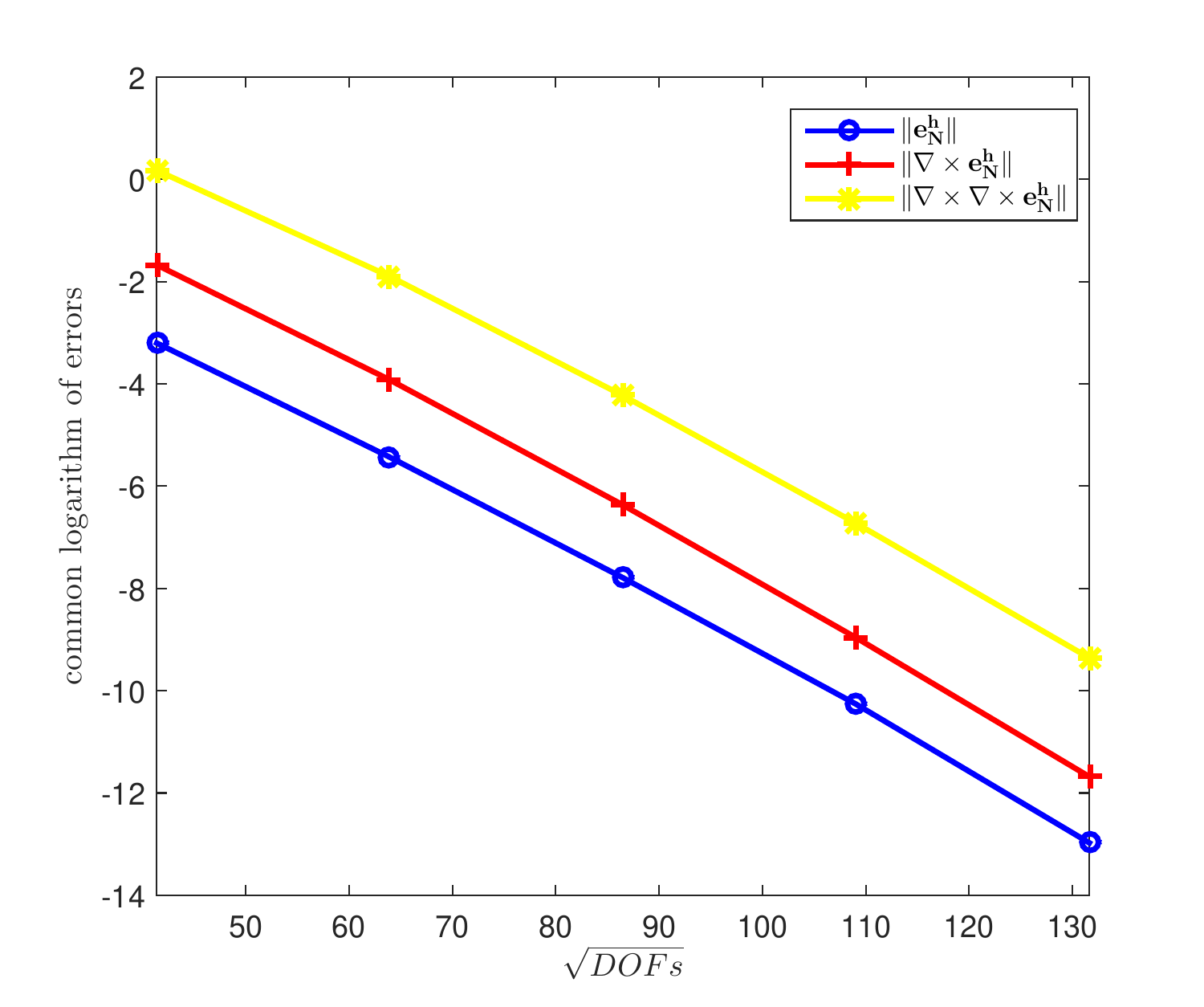}}

  \caption{The $\bm L^2$-errors of $\bm u, \nabla\times\bm u$ and $(\nabla\times)^2\bm u$ versus $\sqrt{DOFs}$ in semi-log scale by  using the   $H(\tc^2)$-conforming quadrilateral spectral elements. }\label{spectralaccuracy}
\end{figure}

To test the $p$-convergence of our method, we set   $L=M=N$ and let  (a)  $h=1$ for one element,  
(b) $h=1/2$  for 4 elements, 
 (c)  $h=1/8$  for  64 elements.   Numerical
 errors  versus  various  DOFs are  shown in Figure \ref{spectralaccuracy} in a semi-logrithm  scale.
Three plots in Figure \ref{spectralaccuracy}  demonstrate the exponential
 orders of convergence of our $H(\curl^2)$-conforming  quadrilateral  spectral element method.
It is noted that $h=1$  provides the optimal  convergence rate, and an obvious decrease in  the convergence rate
is observed as the total number of elements increases.  Moreover, for fixed DOFs,  the accuracy  decrease  monotonously as $h$ decreases
and thus $h=1$ offers the highest accuracy.

\subsection{The quad-curl eigenvalue problem}
We propose the quadrilateral spectral element method to solve the
quad-curl eigenvalue problem which just substitutes $\lambda \bm u$ for the right-hand term $\bm f$  of the source problem \eqref{prob22}.

 As we have seen from Subsection 6.1, the selection of $L$ and $M$   affects the convergence rate for $\pmb u$ only. On the other hand, the accuracy of eigenvalue approximation has much to do with accuracy of computed $\nabla \times \nabla \times \pmb u$. Hence, in this subsection, we always choose $L=M=N$.

The approximation scheme for the  eigenvalue problem  is to find $\lambda_N^h\in \mathbb{R}$ and $(\bm u_N^h, p_N^h)\in \OO{W}_N^{h}\times \OO{S}_N^{h}$, s.t.
\begin{equation}\label{eig-01}
\begin{cases}
\displaystyle a(\bm u_N^h,\bm v_N^h) + b(\bm v_N^h,p_N^h)=\lambda_N^h(\bm u_N^h, \bm v_N^h), & \forall \bm v_N^h\in \OO{W}_N^{h},\\
\displaystyle b(\bm u_N^h,q_N^h)=0, &\forall q_N^h\in \OO{S}_N^{h}.
\end{cases}
\end{equation}
\subsubsection{Square domain}

\begin{figure}[htbp]
  \subfigure[$N=1$;]{\includegraphics[width=0.32\textwidth,height=0.28\textwidth]{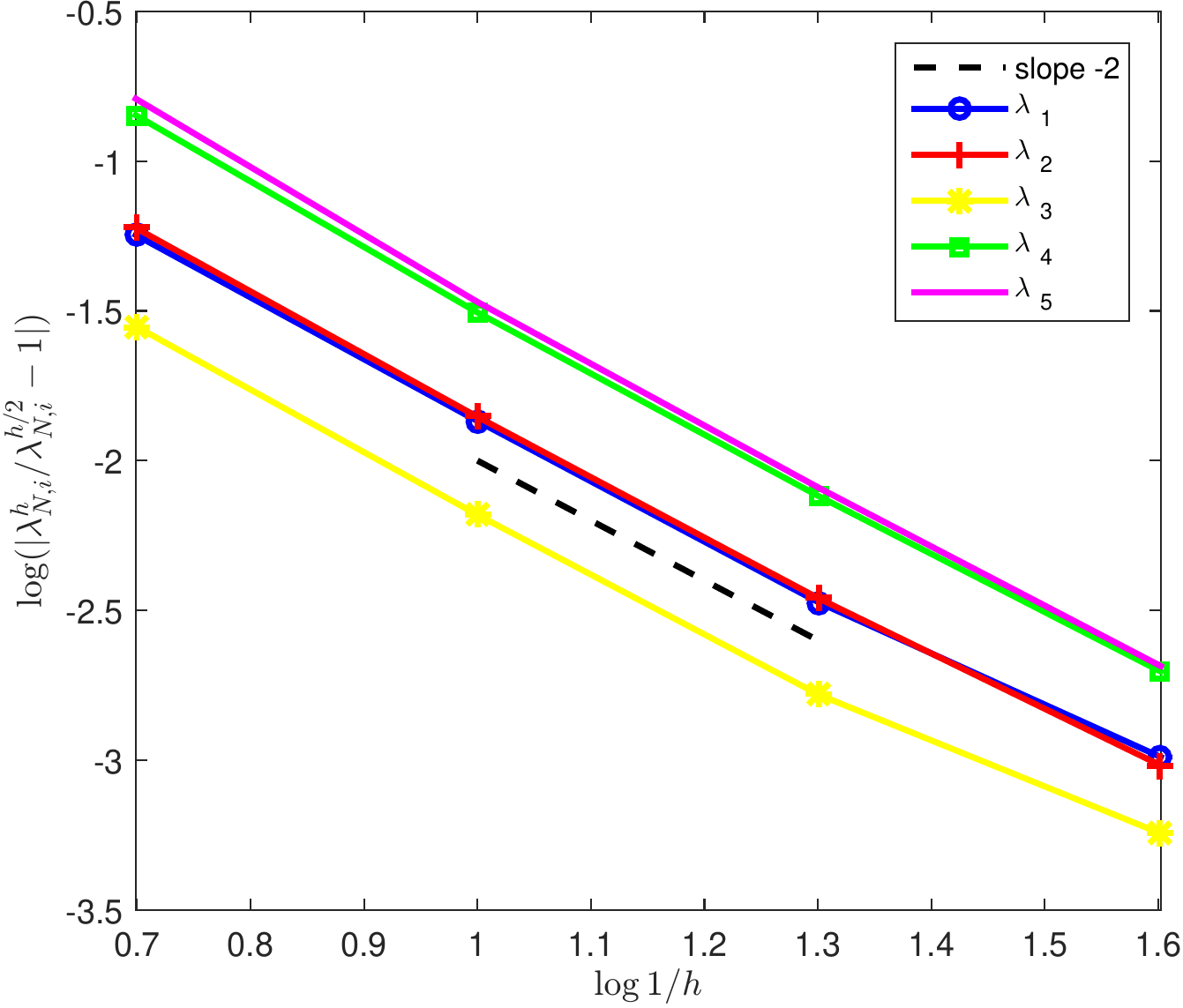}}%
  \subfigure[$N=3$;]{\includegraphics[width=0.32\textwidth,height=0.28\textwidth]{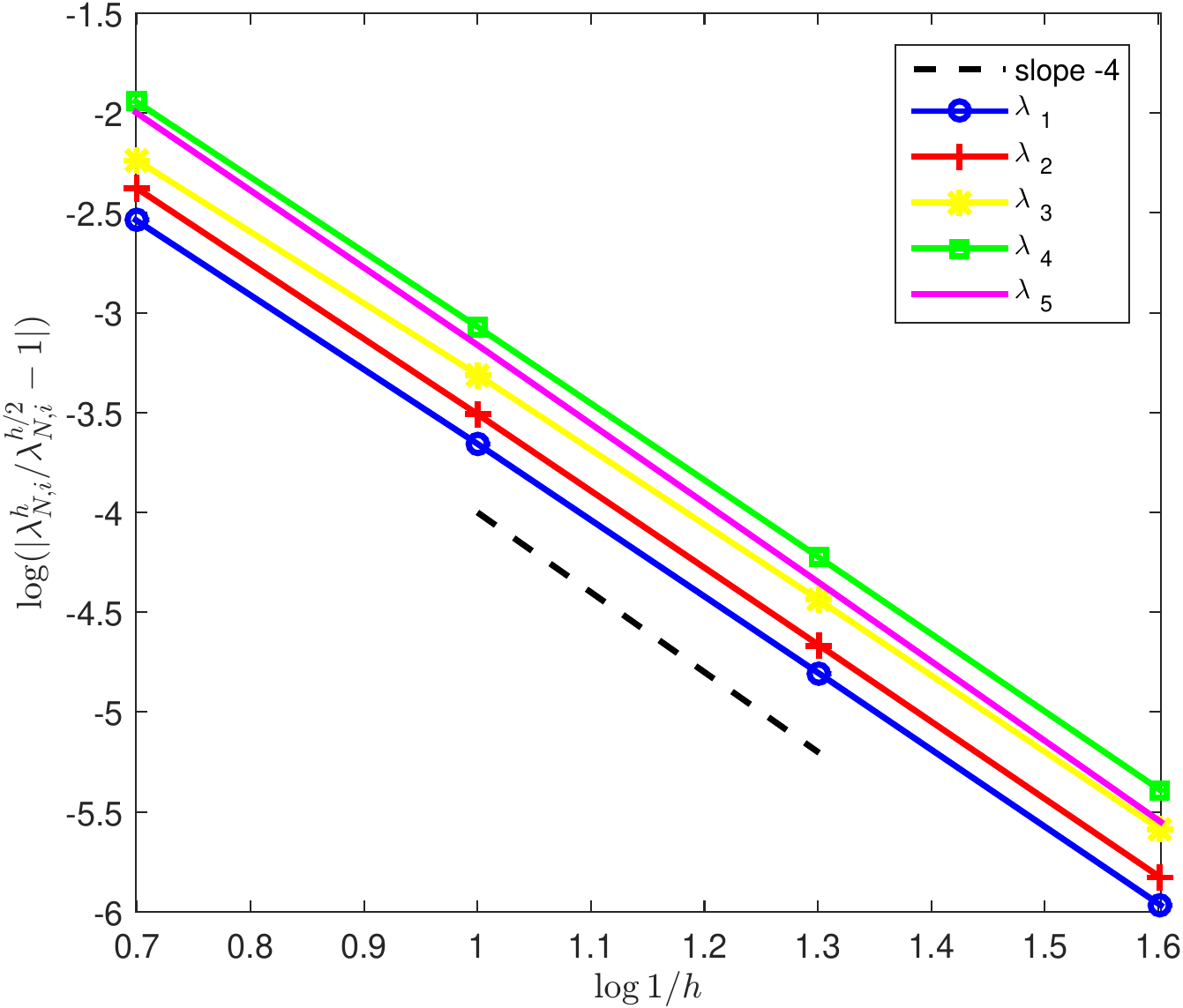}}%
  \subfigure[$N=4$.]{\includegraphics[width=0.32\textwidth,height=0.28\textwidth]{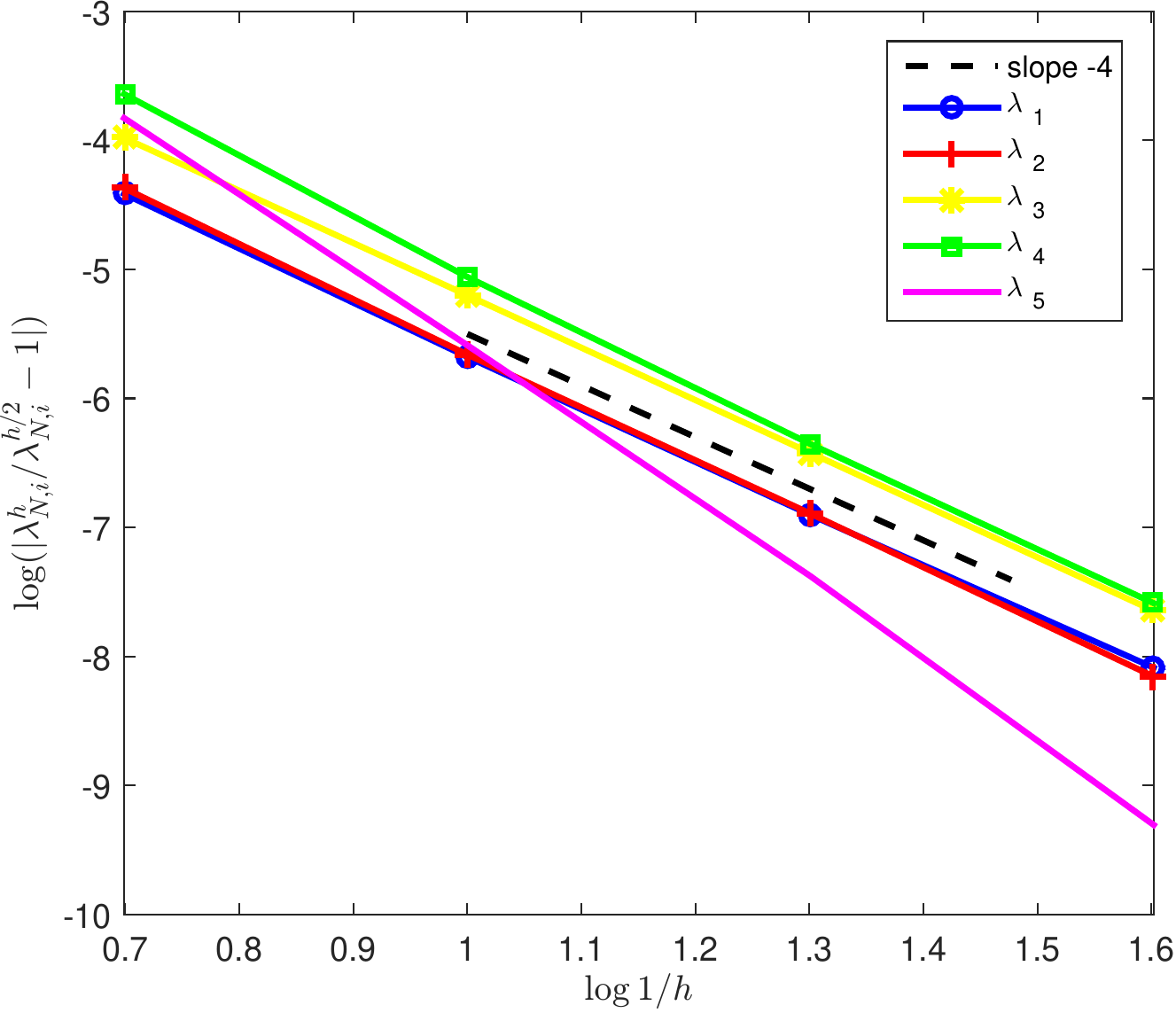}}

  \caption{Eigenvalue errors versus  $h$ for the first five eigenvalues on $[0,1]^2$.}\label{eigfig1}
\end{figure}

 Denote by $\lambda_{i}$ the $i$-th  exact eigenvalue and by $\lambda_{N,i}^h$  the $i$-th numerical eigenvalue. We demonstrate convergence orders of the  first five discrete eigenvalues with various mesh sizes by using
the  simplest  $H(\curl^2)$-conforming  quadrilateral elements in  Table \ref{eigtablower51}.
Relative errors are  then plotted in Figure \ref{eigfig1} (a).
 Second-order $h$-convergence can be observed  in both  Table \ref{eigtablower51} and  Figure \ref{eigfig1} (a)   for each discrete eigenvalue. Noting that an eigenvalue  converges twice as fast as its eigenfunction,  we confirm that
 the convergence  orders of the simplest $H(\curl^2)$-conforming quadrilateral  element method for  the eigenvalue problem
  are consistent with  those for  the source problem reported in  the last subsection.

\begin{table}[h!]
  \centering
  \caption{The first 5  quad-curl  eigenvalues on $[0,1]^2$  and their convergence orders    by using  $H(\tc^2)$-conforming quadrilateral elements with $N=1$.}\label{eigtablower51}
  \begin{tabular}{cccccccccccc}
   \hline
   $h$&$\lambda_{1,1}^h$ & order&$\lambda_{1,2}^h$&  order& $\lambda_{1,3}^h$&order&$\lambda_{1,4}^h$ & order&$\lambda_{1,5}^h$&order&\\
   \hline
   $1/5$&761.635&-&764.313&-&2437.229&-&5063.90&-&6097.93&-\\
   $1/10$&720.783&2.06&721.198&2.09&2370.70&2.03&4433.19&2.16&5250.20&2.25\\
   $1/20$&711.142&2.01&711.241&2.02&2355.17&1.99&4298.94&2.04&5078.26&2.06\\
   $1/40$&708.763&2.00&708.787&2.00&2351.29&1.99&4266.53&2.01&5037.42&2.01\\
      $1/80$&708.169&-&708.175&-&2350.31&-&4258.49&-&5027.34&-\\
      \hline
  \end{tabular}
 \end{table}

\begin{table}[h]
	\centering
	\caption{The first 5  quad-curl  eigenvalues on $[0,1]^2$  and their convergence orders   by using the $H(\tc^2)$-conforming elements with $N=3$. } \label{tab3}
	\begin{tabular}{ccccccccccc}
		\hline
		$h$    &$\lambda^h_{3,1}$&order&$\lambda^h_{3,2}$&order&$\lambda^h_{3,3}$&order&$\lambda^h_{3,4}$&order&$\lambda^h_{3,5}$&order\\
		\hline
		$1/5$&710.200&-&711.203&-&2364.96&-&4308.39&-&5078.35&-\\
		$1/10$&708.139&3.72&708.207&3.77&2351.22&3.58&4259.70&3.75&5027.69&3.87\\
	    $1/20$&707.983&3.82&707.987&3.84&2350.07&3.75&4256.08&3.84&5024.23&3.95\\
$1/40$&707.972&3.86&707.972&3.87&2349.99&3.82&4255.83&3.88&5024.00&3.98\\
	    $1/80$&707.971&-&707.971&-&2349.98&-&4255.81&-&5023.99&-\\
		\hline
	\end{tabular}
\end{table}

{\addtolength{\tabcolsep}{-1pt}
\begin{table}[h]
	\centering
	\caption{The first 5  quad-curl  eigenvalues on $[0,1]^2$  and their convergence orders   by using the $H(\tc^2)$-conforming elements with $N=4$. } \label{tab4}
	\begin{tabular}{ccccccccccccc}
		\hline
		$h$    &$\lambda^h_{4,1}$&order&$\lambda^h_{4,2}$&order&$\lambda^h_{4,3}$&order&$\lambda^h_{4,4}$&order&$\lambda^h_{4,5}$&order\\
		\hline
		$1/5$&708.0004&-&708.0034&-&2350.2475&-&4256.8267&-&5024.7537&-\\
		$1/10$&707.9731&4.20&707.9732&4.29&2350.0016&4.05&4255.8534&4.71&5024.0055&5.85\\
	    $1/20$&707.9716&4.07&707.9716&4.10&2349.9868&4.06&4255.8162&4.30&5023.9926&5.95\\
$1/40$&707.9715&3.93&707.9715&4.19&2349.9859&4.06&4255.8143&4.09&5023.9923&6.40\\
	    $1/80$&707.9715&-&707.9715&-&2349.9859&-&4255.8142&-&5023.9923&-\\
		\hline
	\end{tabular}
\end{table}
}

 Further,  let us set $h=1/5$ as the initial mesh size and carry out  numerical experiments on the $h$-version of the $H(\curl^2)$-conforming quadrilateral elements for   \eqref{eig-01}  with $N=3,4$. The first five discrete eigenvalues  are shown in
 Table \ref{tab3} for $N=3$ and Table \ref{tab4} for $N=4$, and  their relative errors are depicted in  Figure \ref{eigfig1} (b) and  (c).
 From these tables and figures, one readily finds that all  five  discrete eigenvalues converge at the  full order around $4$
 for $N=3$. At the same time, only the fifth discrete eigenvalue converges at the full order  around $6$,  while  the first four eigenfunctions
 have still a convergence order  slightly larger than $4$ owing to the limited regularity of their eigenfunctions.

\begin{figure}[htbp]
  \subfigure[$h=1$;]{\includegraphics[width=0.32\textwidth]{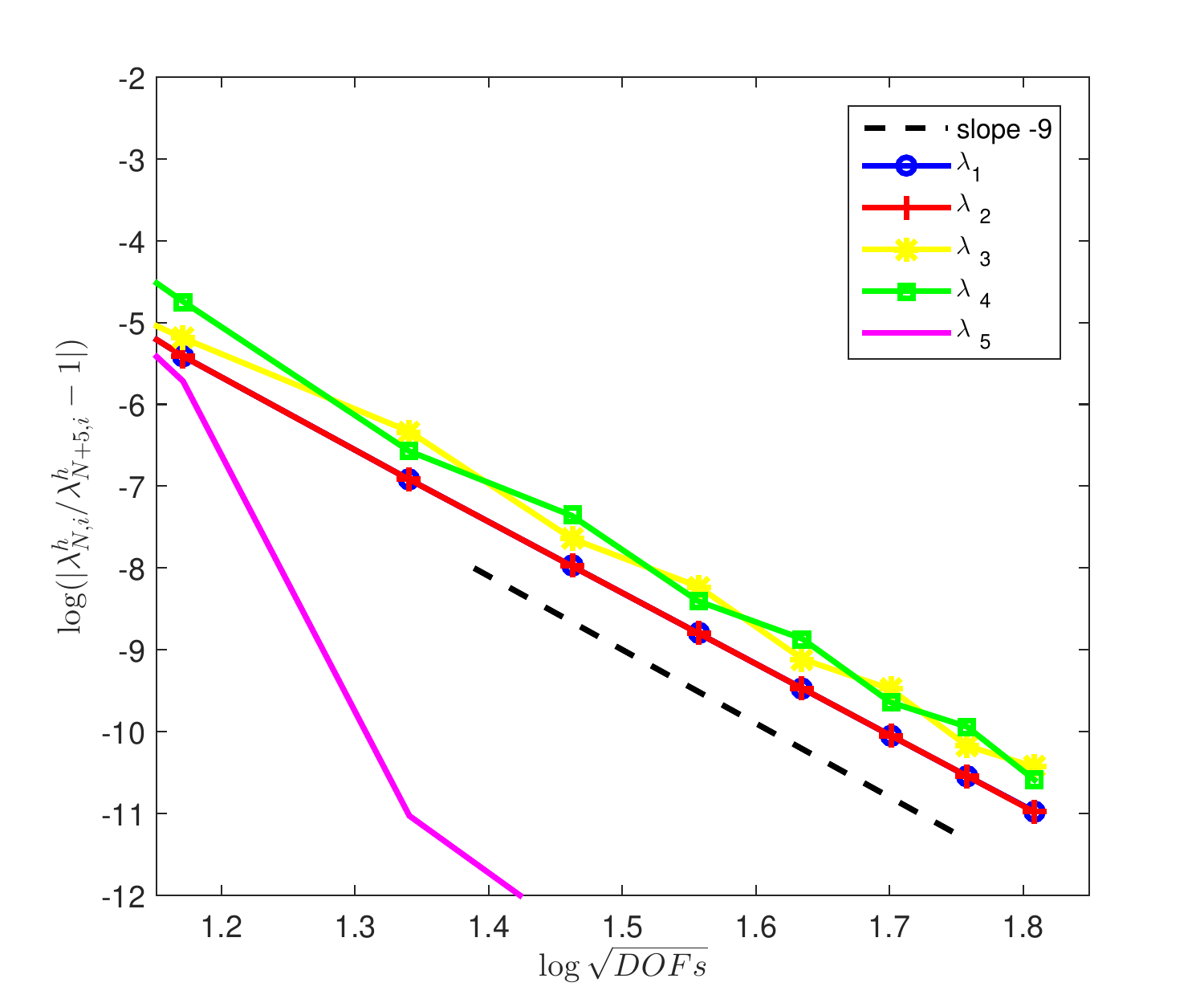}}%
   \subfigure[$h=1/2$;]{\includegraphics[width=0.32\textwidth]{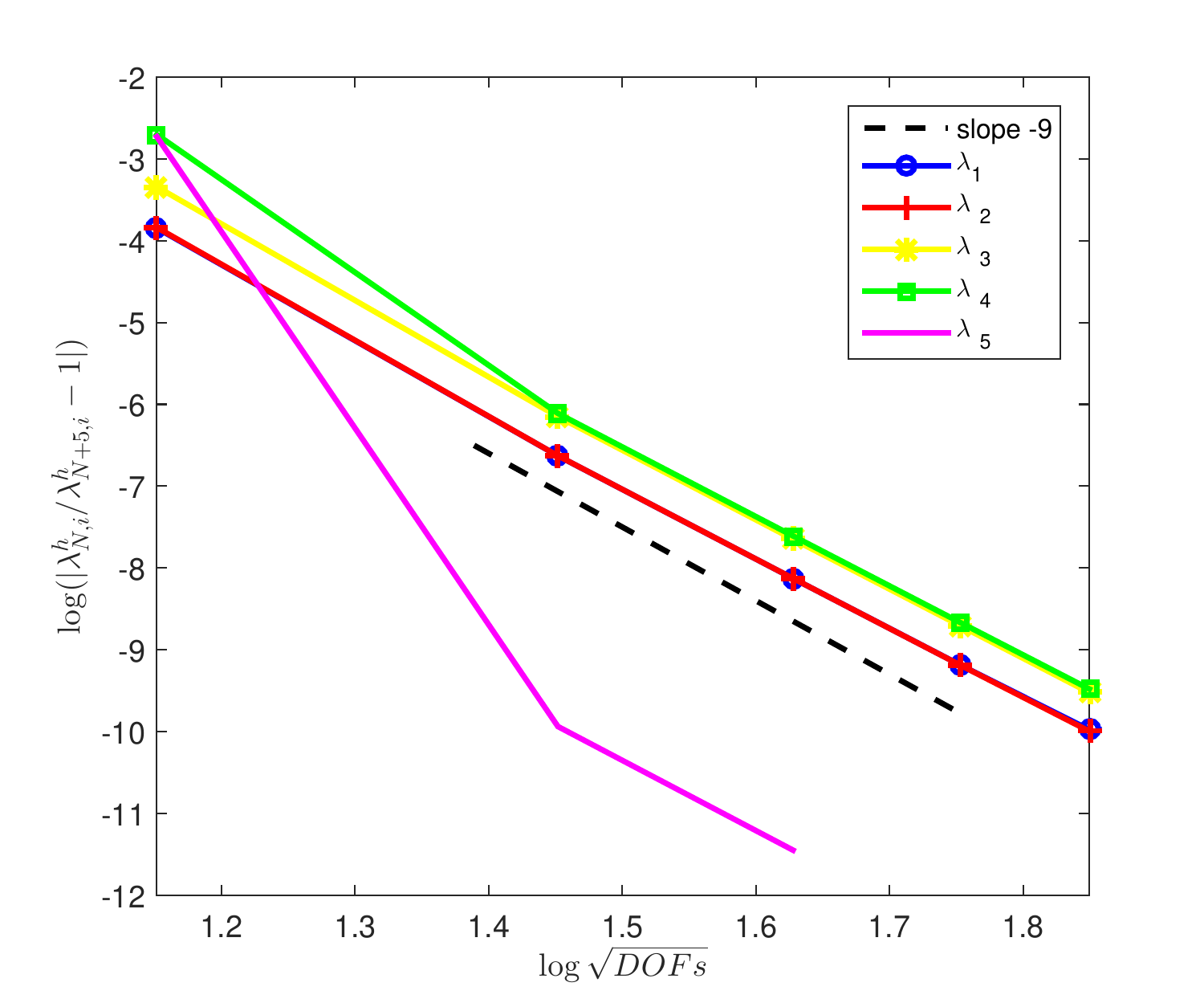}}%
    \subfigure[$h=1/10$.]{\includegraphics[width=0.32\textwidth]{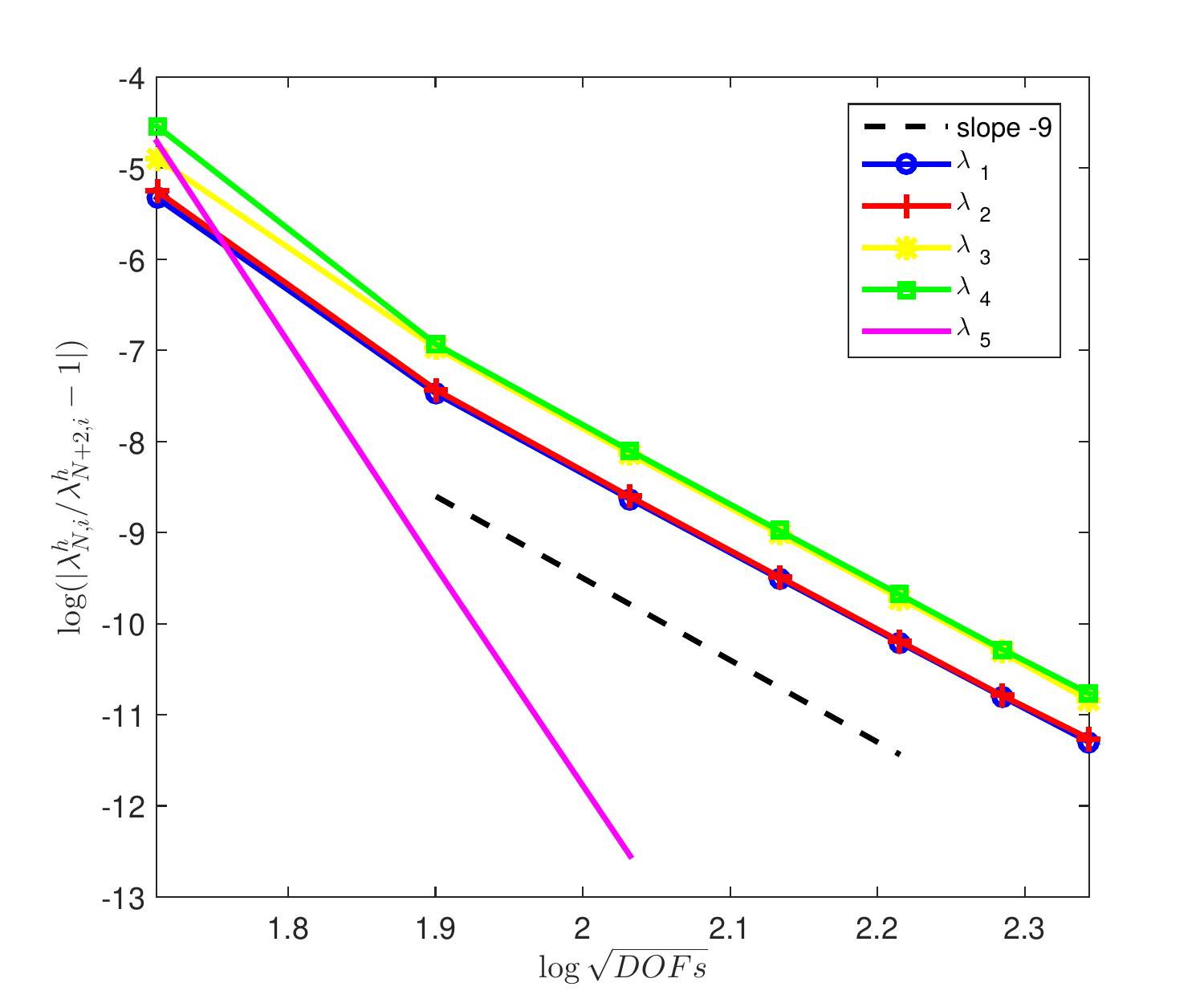}}
  \caption{Eigenvalue errors versus  $\sqrt{DOFs}$ of the $H(\tc^2)$-conforming quadrilateral spectral  method on $[0,1]^2$.}\label{spectral}
\end{figure}

 Next, we let  (a)  $h=1$ for one element,
(b) $h=1/2$  for 4 elements,
 (c)  $h=1/10$  for  100 elements,  then  examine the $p$-convergence of  our  $H(\curl^2)$-conforming quadrilateral elements.
 Errors of numerical eigenvalues versus  various  DOFs are plotted  in Figure \ref{spectral} in  log-log  scale.
Distinct to those for infinitely smooth problems,  our $H(\curl^2)$-conforming quadrilateral spectral element methods
for quad-curl eigenvalue problems  only have  algebraic orders of convergence.
  Indeed,
   the quad-curl eigenvalue problem is essentially  a sixth order partial differential equation,  and
singularities  shall occur even on a square domain, so that only limited convergence orders can be obtained in both the $p$- and $h$-versions  of  our method.
Nevertheless,
the convergence rates for a fixed eigenvalue are independent of the total number of elements used in our $p$-version spectral element methods. Indeed, they are  twice as   high as those of the $h$-version for $N\ge 4$.

\subsubsection{L-shaped domain}
The quad-curl eigenvalue problem   on an L-shaped domain $\Omega=(0,1)\times(0,1)/[0.5,1)\times[0.5,1)$ is also considered.  Due to the strong singularity of the domain, convergence
rate for the first eigenvalue deteriorates to around $h^{4/3}$ in the $h$-version (see Table \ref{tab-LL-01}).
While, it is observed that  the convergence rate in the $p$-version  with  $h=1/6$ is nearly $N^{-3.5}$   (see Figure \ref{L-shape-01}).
Once again, these  reflect the correctness and efficiency of our  $H(\tc^2)$-conforming quadrilateral
spectral element method.

\vspace*{2em}
\begin{minipage}{0.6\textwidth}
\centering
  \captionsetup{width=1\textwidth}
\captionof{table}{The first quad-curl eigenvalue  and the convergence order  by  $H(\curl^2)$-conforming quadrilateral spectral elements with $N=4$ on the L-shaped domain.}
\begin{tabular}{cccc}
  \hline
  $h$ & $\lambda^h_{4,1}$ & error&order \\
  \hline
  $1/4$ &534.46527767676&8.920012568e-04&-\\
  $1/8$ &534.94202137614&4.411857535e-04&1.0156\\
  $1/16$ &535.17803017493&1.814911829e-04&1.2815\\
  $1/32$ &535.27516026869&7.262977955e-05&1.3213\\
  $1/64$ &535.31403718558&2.889401638e-05&1.3298\\
  $1/128$ &535.32950455814&-&-\\
  \hline
\end{tabular}\label{tab-LL-01}
\end{minipage}\hfill%
\begin{minipage}{0.4\textwidth}

  \centering
  \includegraphics[width=0.9\textwidth,height=0.7\textwidth]{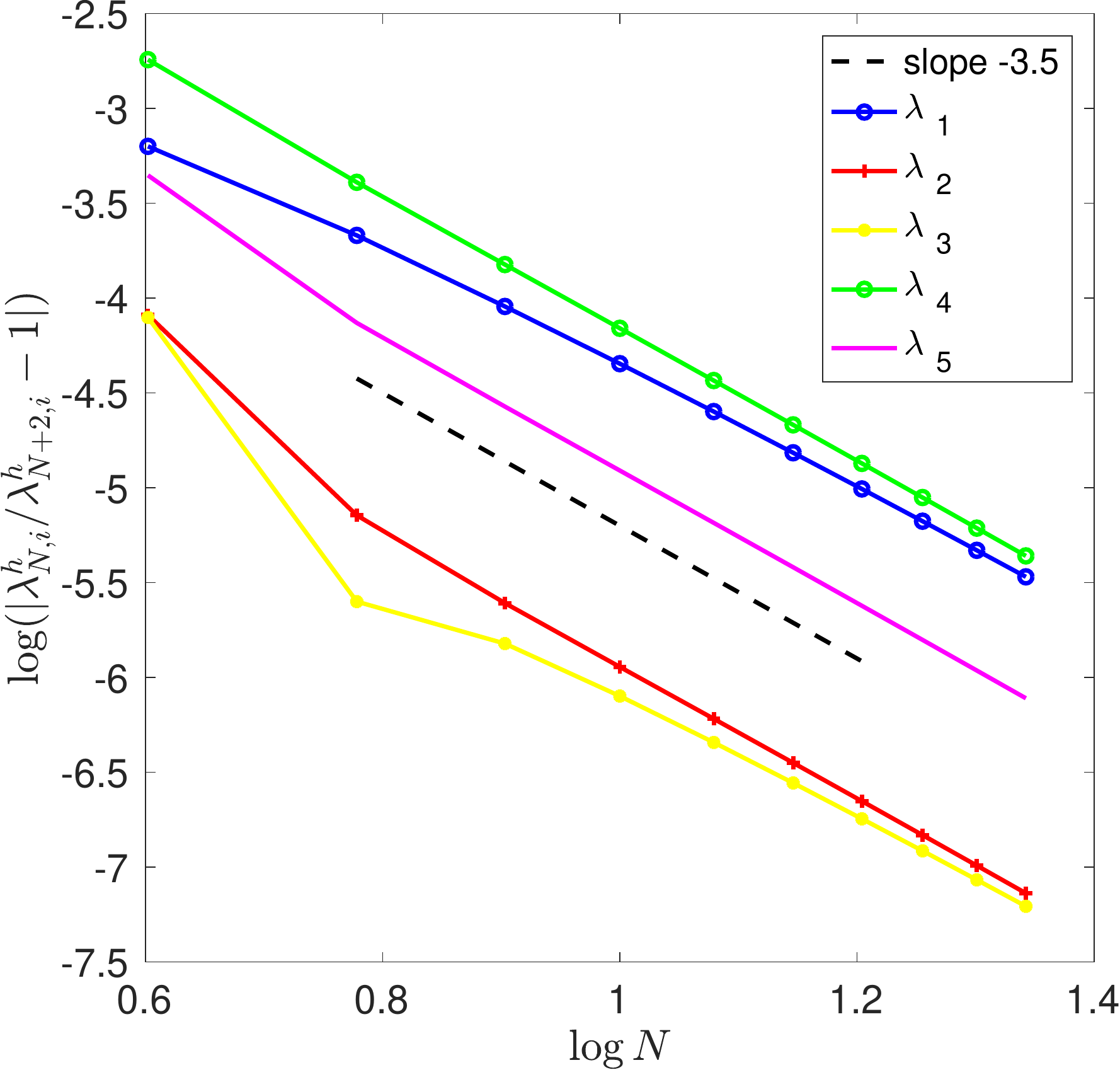}
  \captionsetup{width=0.9\textwidth}
   \captionof{figure}{Eigenvalue errors versus $N$ of the $H(\curl^2)$-conforming quadrilateral spectral element method on the L-shape domain.}\label{L-shape-01}
\end{minipage}\hspace*{\fill}
\vspace*{1em}

\

\section{Conclusion}
 In this paper, we have constructed
$H(\tc^2)$-conforming basis  functions  on an arbitrary  convex quadrilateral  element using
the bilinear mapping from the reference square onto the physical element   together with the contravariant transformation of vector fileds.
These hierarchical basis functions are explicitly formulated in  generalized Jacobi polynomials with the indices of $(-1,-1) $ and $(-2, -2)$,   which are  easier to generalize to higher order approximation scheme than  the  Lagrangian bases \cite{Zhang2019H}. Numerical results show that our $H(\tc^2)$-conforming spectral elements are efficient   and  can achieve an exponential order of $p$-convergence
 and an optimal order,  $\mathcal{O}(h^{N-1})$,   of  $h$-convergence   for  source problems  measured in the $H(\tc^2;\Omega)$-norm.
However,  owing to the strong singularities,   our $H(\tc^2)$-conforming spectral elements  get only  limited orders of convergence
in both the $p$- and $h$-versions  for eigenvalue problems.

\begin{appendix}

\section{Proof of Lemma \ref{lowerorderpolynomials}}
\label{AppendixA}
$Proof\ of\ Lemma\ \emph{\ref{lowerorderpolynomials}}.$
Based on \eqref{Qfunc1}, \eqref{Qfunc2}, \eqref{Qfunc3}, and, \eqref{Qfunc4}, one can obtain
\begin{align*}
  &-x_{14}\hat {\bm \phi}_{0,0}+x_{21}\hat {\bm \phi}_{1,0}+x_{32}\hat {\bm \phi}_{1,1}-x_{43}\hat {\bm \phi}_{0,1}\\
 =&\left(\frac{x_{21}+x_{34}}{4}+\frac{-x_{21}+x_{34}}{4}\hat y,\frac{x_{41}+x_{32}}{4}+\frac{-x_{41}+x_{32}}{4}\hat x\right)^\tr\\
 =&B_K^\tr(1,0)^\tr,
\end{align*}
which implies \eqref{low-01}. Similarly,  \eqref{low-02} can be yielded.

Furthermore, we also have
\begin{align}\label{lemm-pro-01}
\begin{split}
&-\frac{(x_1^2-x_4^2)}{2}\hat {\bm \phi}_{0,0}+\frac{(x_2^2-x_1^2)}{2}\hat {\bm \phi}_{1,0}+\frac{(x_3^2-x_2^2)}{2}\hat {\bm \phi}_{1,1}-\frac{(x_4^2-x_3^2)}{2}\hat {\bm \phi}_{0,1}\\
=&\Big(\frac{(x_1^2-x_2^2+x_3^2-x_4^2)\hat y}{8} +\frac{- x_1^2+x_2^2+ x_3^2- x_4^2}{8},\\
&\frac{(x_1^2-x_2^2+x_3^2-x_4^2)\hat x}{8} +\frac{- x_1^2-x_2^2+ x_3^2+ x_4^2}{8}\Big)^\tr.
\end{split}
\end{align}
By \eqref{Qfunc1}-\eqref{Qfunc4}, one finds
\begin{align}\label{lemm-pro-02}
\begin{split}
&\frac{x_{14}^2}{2}\hat {\bm \phi}_{0,2}+\frac{x_{21}^2}{2}\hat {\bm \phi}_{2,0}+\frac{x_{32}^2}{2}\hat {\bm \phi}_{1,2}+\frac{x_{43}^2}{2}\hat {\bm \phi}_{2,1}\\
=&\Bigg(-\frac{x_{41}^2 (\hat y^2-1)}{16}+\frac{x_{21}^2 \hat x (1-\hat y)}{8}+\frac{ x_{32}^2 ( \hat y^2-1)}{16}+\frac{ x_{43}^2 \hat x (1+\hat y)}{8},\\
&\frac{x_{41}^2 (1- \hat x) \hat y}{8}-\frac{x_{21}^2 ( \hat x^2-1)}{16}+\frac{x_{32}^2 (1+ \hat x) \hat y}{8}+\frac{ x_{43}^2 ( \hat x^2-1)}{16}\Bigg)^\tr,
\end{split}
\end{align}
and
\begin{align}\label{lemm-pro-03}
\begin{split}
&\frac{(x_{32}+x_{14})^2}{2}\hat {\bm \phi}_{2,2} \\
=&\Bigg(\frac{(x_{32}+x_{14})^2\hat x(\hat y^2-1)}{16},\frac{(x_{32}+x_{14})^2(\hat x^2-1)\hat y}{16}\Bigg)^\tr.
\end{split}
\end{align}
Adding \eqref{lemm-pro-01}-\eqref{lemm-pro-03} up, then
\begin{align*}
&  -\frac{(x_1^2-x_4^2)}{2}\hat {\bm \phi}_{0,0}+\frac{(x_2^2-x_1^2)}{2}\hat {\bm \phi}_{1,0}+\frac{(x_3^2-x_2^2)}{2}\hat {\bm \phi}_{1,1}-\frac{(x_4^2-x_3^2)}{2}\hat {\bm \phi}_{0,1} + \frac{x_{14}^2}{2}\hat {\bm \phi}_{0,2}\nonumber
\\
&+\frac{x_{21}^2}{2}\hat {\bm \phi}_{2,0}+\frac{x_{32}^2}{2}\hat {\bm \phi}_{1,2}+\frac{x_{43}^2}{2}\hat {\bm \phi}_{2,1}+\frac{(x_{32}+x_{14})^2}{2}\hat {\bm \phi}_{2,2} \\
=&\Bigg(\left(\frac{x_{21}}{2}\frac{1-\hat y}{2}+\frac{x_{34}}{2}\frac{1+\hat y}{2}\right)\left(\sigma_1(\hat x,\hat y)x_1+\sigma_2(\hat x,\hat y)x_2+\sigma_3(\hat x,\hat y)x_3+\sigma_4(\hat x,\hat y)x_4\right),\\
&\left(\frac{x_{41}}{2}\frac{1-\hat x}{2}+\frac{x_{32}}{2}\frac{1+\hat x}{2}\right)\left(\sigma_1(\hat x,\hat y)x_1+\sigma_2(\hat x,\hat y)x_2+\sigma_3(\hat x,\hat y)x_3+\sigma_4(\hat x,\hat y)x_4\right)\Bigg)^\tr\\
=&B_K^\tr(x,0)^\tr.
\end{align*}
Hence, we arrive at \eqref{low-03}. If analyzing a bit, we can find \eqref{low-05}.

Now, let us prove \eqref{low-04} and \eqref{low-06}. Some equations are calculated as follows:
\begin{align}
& -\frac{(x_1+x_4)y_{14}}{2}\hat {\bm \phi}_{0,0}+\frac{(x_2+x_1)y_{21}}{2}\hat {\bm \phi}_{1,0}+\frac{(x_3+x_2)y_{32}}{2}\hat {\bm \phi}_{1,1}-\frac{(x_4+x_3)y_{43}}{2}\hat {\bm \phi}_{0,1}\nonumber\\
=&\Bigg(\frac{y_{42}x_{13}-y_{13}x_{42}}{64}\hat y(\hat y^2-1)(3\hat x^2-5)+\frac{(x_2+x_1) y_{21} }{8}(1-\hat y)-\frac{(x_4+x_3)y_{43}}{8}(1+\hat y),\nonumber\\
&\frac{y_{24}x_{13}-y_{13}x_{24}}{64}\hat x(\hat x^2-1)(3\hat y^2-5)+\frac{(x_4+x_1) y_{41} }{8}(1-\hat x)-\frac{(x_2+x_3)y_{23}}{8}(1+\hat x)\Bigg)^\tr,\label{lemma-pro-11}
\end{align}
\begin{align}
&\frac{x_{14}y_{14}}{2}\hat {\bm \phi}_{0,2}+\frac{x_{21}y_{21}}{2}\hat {\bm \phi}_{2,0}+\frac{x_{32}y_{32}}{2}\hat {\bm \phi}_{1,2}+\frac{x_{43}y_{43}}{2}\hat {\bm \phi}_{2,1}\nonumber\\
=&\Bigg(-\frac{x_{41} y_{41} (\hat y^2-1)}{16}+\frac{ x_{21} y_{21} \hat x (1- \hat y)}{8}+\frac{x_{32} y_{32} (\hat y^2-1)}{16}+\frac{x_{43} y_{43} \hat x (1+ \hat y)}{8},\nonumber\\
 &\frac{x_{41} y_{41} (1- \hat x) \hat y}{8}-\frac{ x_{21} y_{21} ( \hat x^2-1)}{16}+\frac{ x_{32} y_{32} (1+ \hat x) \hat y}{8}+\frac{x_{43} y_{43} ( \hat x^2-1)}{16}\Bigg)^\tr,\label{lemma-pro-21}
\end{align}
and,
\begin{align}
&\frac{(x_{32}+x_{14})(y_{32}+y_{14})}{2}\hat {\bm \phi}_{2,2}\nonumber\\
=&\Bigg(\frac{ (x_{32}+x_{14}) (y_{32}+y_{14}) }{4}\frac{ \hat y^2-1}{4}\hat x,\frac{(x_{32}+x_{14}) (y_{32}+y_{14}) }{4}\frac{\hat x^2-1}{4} \hat y\Bigg)^\tr.\label{lemma-pro-31}
\end{align}
Replacing  $l_i, s_i$ with $x_i, y_i, i=1,\cdots, 4$, one can arrive at
\begin{align}
&\hat {\bm \psi}_{0,0}+\hat {\bm \psi}_{0,1}+\hat {\bm \psi}_{1,0}+\hat {\bm \psi}_{1,1}\nonumber\\
=&\Bigg(\frac{3(1-\hat y)(1+\hat y)}{64}\Big((y_{42}x_{13}-y_{13}x_{42})(\hat x^2-\frac{5}{3})\hat y+\frac{4}{3}(y_{32}x_{14}-y_{14}x_{32})\Big),\nonumber\\
&\frac{3(\hat x-1)(1+\hat x)}{64}\Big((y_{42}x_{13}-y_{13}x_{42})(\hat y^2-\frac{5}{3})\hat x+\frac{4}{3}(y_{43}x_{12}-y_{12}x_{43})\Big)\Bigg)^\tr.\label{lemma-pro-41}
\end{align}
Summing \eqref{lemma-pro-11}-\eqref{lemma-pro-41}, it gives
\begin{align*}
 & -\frac{(x_1+x_4)y_{14}}{2}\hat {\bm \phi}_{0,0}+\frac{(x_2+x_1)y_{21}}{2}\hat {\bm \phi}_{1,0}+\frac{(x_3+x_2)y_{32}}{2}\hat {\bm \phi}_{1,1}-\frac{(x_4+x_3)y_{43}}{2}\hat {\bm \phi}_{0,1}
\nonumber\\
&+ \frac{x_{14}y_{14}}{2}\hat {\bm \phi}_{0,2}+\frac{x_{21}y_{21}}{2}\hat {\bm \phi}_{2,0}+\frac{x_{32}y_{32}}{2}\hat {\bm \phi}_{1,2}+\frac{x_{43}y_{43}}{2}\hat {\bm \phi}_{2,1}
\nonumber\\
&+\frac{(x_{32}+x_{14})(y_{32}+y_{14})}{2}\hat {\bm \phi}_{2,2}+(\hat {\bm \psi}_{0,0}+\hat {\bm \psi}_{0,1}+\hat {\bm \psi}_{1,0}+\hat {\bm \psi}_{1,1})\\
=&\Bigg(\left(\frac{y_{21}}{2}\frac{1-\hat y}{2}+\frac{y_{34}}{2}\frac{1+\hat y}{2}\right)\left(\sigma_1(\hat x,\hat y)x_1+\sigma_2(\hat x,\hat y)x_2+\sigma_3(\hat x,\hat y)x_3+\sigma_4(\hat x,\hat y)x_4\right),\\
&\left(\frac{y_{41}}{2}\frac{1-\hat x}{2}+\frac{y_{32}}{2}\frac{1+\hat x}{2}\right)\left(\sigma_1(\hat x,\hat y)x_1+\sigma_2(\hat x,\hat y)x_2+\sigma_3(\hat x,\hat y)x_3+\sigma_4(\hat x,\hat y)x_4\right)\Bigg)^\tr\\
=&B_K^\tr(0,x)^\tr,
\end{align*}
which states \eqref{low-04}. Similarly, we derive \eqref{low-06}.
\eqref{low-07} can be deduced from \eqref{low-04} since $\hat \nabla \times \hat{\bm {\bm \phi}}_{m,n}=0$, when $(m,n)\neq(i,j), \{i,j\}\in\{0,1\}$. Now we finish the proof.

\end{appendix}	
	
\bibliographystyle{plain}
\bibliography{reference}

\begin{thebibliography}{10}

\bibitem{Arnold2018Finite}
D.~Arnold.
\newblock {\em Finite Element Exterior Calculus}.
\newblock Society for Industrial and Applied Mathematics, Philadelphia, PA,
  2018.

\bibitem{Belgacem1999Spectral}
F.~Ben Belgacem and C.~Bernardi.
\newblock Spectral element discretization of the {M}axwell equations.
\newblock {\em Mathematics of Computation}, 68(228):1497--1520, 1999.

\bibitem{Brenner2019Multigrid}
S.~C. Brenner, J.~Cui, and L.~Sung.
\newblock Multigrid methods based on {H}odge decomposition for a quad-curl
  problem.
\newblock {\em Computational Methods in Applied Mathematics}, 19(2):215--232,
  2019.

\bibitem{Brenner2017Hodge}
S.~C. Brenner, J.~Sun, and L.~Sung.
\newblock Hodge decomposition methods for a quad-curl problem on planar
  domains.
\newblock {\em Journal of Scientific Computing}, 73(2-3):495--513, 2017.

\bibitem{Cai2013}
W.~Cai.
\newblock {\em Computational Methods for Electromagnetic Phenomena}.
\newblock Cambridge University Press, New York, 2013.

\bibitem{Cakoni2017A}
F.~Cakoni and H.~Haddar.
\newblock A variational approach for the solution of the electromagnetic
  interior transmission problem for anisotropic media.
\newblock {\em Inverse Problems and Imaging}, 1(3):443--456, 2017.

\bibitem{Cohen2002Higher}
G.~Cohen.
\newblock {\em Higher-Order numerical methods for transient wave equations}.
\newblock Springer-Verla, New York, 2002.

\bibitem{Guo2009Generalized}
B.~Guo, J.~Shen, and L.~Wang.
\newblock Generalized jacobi polynomials/functions and their applications.
\newblock {\em Applied Numerical Mathematics}, 59(5):1011--1028, 2009.

\bibitem{Guo2010Composite}
B.~Guo and T.~Wang.
\newblock Composite {L}aguerre-{L}egendre spectral method for fourth-order
  exterior problems.
\newblock {\em Journal of Scientific Computing}, 44(3):255--285, 2010.

\bibitem{Qingguo2012A}
Q.~Hong, J.~Hu, S.~Shu, and J.~Xu.
\newblock A discontinuous {G}alerkin method for the fourth-order curl problem.
\newblock {\em Journal of Computational Mathematics}, 30(6):565--578, 2012.

\bibitem{Hu2020simple}
K~Hu, Q.~Zhang, and Z.~Zhang.
\newblock Simple curl-curl-conforming finite elements in two dimensions.
\newblock {\em arXiv}, 2020.

\bibitem{Li2019C1}
H.~Li, W.~Shan, and Z.~Zhang.
\newblock ${C}^1$-conforming quadrilateral spectral element method for
  fourth-order equations.
\newblock {\em Communications on Applied Mathematics and Computation},
  1(3):403--434, 2019.

\bibitem{Liu2010A}
Y.~Liu, J.~Lee, X.~Tian, and Q.~Liu.
\newblock A spectral-element time-domain solution of {M}axwell's equations.
\newblock {\em Microwave \& Optical Technology Letters}, 48(4):673--680, 2006.

\bibitem{Monk2012Finite}
P.~Monk and J.~Sun.
\newblock Finite element methods for {M}axwell's transmission eigenvalues.
\newblock {\em SIAM Journal on Scientific Computing}, 34(3):B247--B264, 2012.

\bibitem{Na2015Mixed}
L.~Na, L.~Tob\'on, Y.~Zhao, Y.~Tang, and Q.~Liu.
\newblock Mixed spectral-element method for 3-d {M}axwell's eigenvalue problem.
\newblock {\em IEEE Transactions on Microwave Theory \& Techniques},
  63(2):317--325, 2015.

\bibitem{Nicaise2018Singularities}
S.~Nicaise.
\newblock Singularities of the quad curl problem.
\newblock {\em Journal of Differential Equations}, 264(8):5025--5069, 2018.

\bibitem{Patera1984A}
A.T. Patera.
\newblock A spectral element method for fluid dynamics: Laminar flow in a
  channel expansion.
\newblock {\em Journal of Computational Physics}, 54(3):468--488, 1984.

\bibitem{Shen2009A}
J.~Shen, L.~Wang, and H.~Li.
\newblock A triangular spectral element method using fully tensorial rational
  basis functions.
\newblock {\em SIAM Journal on Numerical Analysis}, 47(3):1619--1650, 2009.

\bibitem{Sun2016A}
J.~Sun.
\newblock A mixed {FEM} for the quad-curl eigenvalue problem.
\newblock {\em Numerische Mathematik}, 132(1):185--200, 2016.

\bibitem{Sun_2013}
J.~Sun and L.~Xu.
\newblock Computation of {M}axwell's transmission eigenvalues and its
  applications in inverse medium problems.
\newblock {\em Inverse Problems}, 29(10):104013, 2013.

\bibitem{Sun2019Awg}
J.~Sun, Q.~Zhang, and Z.~Zhang.
\newblock A curl-conforming weak {G}alerkin method for the quad-curl problem.
\newblock {\em BIT Numerical Mathematics}, 59:1093--1114, 2019.

\bibitem{Sun2016Finite}
J.~Sun and A.~Zhou.
\newblock {\em Finite Element Methods for Eigenvalue Problems}.
\newblock Chapman and Hall/CRC, Boca Raton, FL, 2016.

\bibitem{Sun2018Multigrid}
Z.~Sun, J.~Cui, F.~Gao, and C.~Wang.
\newblock Multigrid methods for a quad-curl problem based on ${C}^0$ interior
  penalty method.
\newblock {\em Computers \& Mathematics with Applications}, 76(9):2192--2211,
  2018.

\bibitem{Szeg1939Orthogonal}
G.~Szeg.
\newblock {\em Orthogonal polynomials}, volume~23.
\newblock American Mathematical Society, 1939.

\bibitem{Wang2019A}
C.~Wang, Z.~Sun, and J.~Cui.
\newblock A new error analysis of a mixed finite element method for the
  quad-curl problem.
\newblock {\em Applied Mathematics and Computation}, 349:23--38, 2019.

\bibitem{WangToappear}
L.~Wang, Q.~Zhang, J.~Sun, and Z.~Zhang.
\newblock A priori and a posteriori error estimates for the quad-curl
  eigenvalue problem.
\newblock {\em arXiv}, 2020.

\bibitem{wangzhangzhang2018}
L.~Wang, Q.~Zhang, and Z.~Zhang.
\newblock Superconvergence analysis and {PPR} recovery of arbitrary order edge
  elements for maxwell's equations.
\newblock {\em Journal of Scientific Computing}, 78:1207--1230, 2019.

\bibitem{Sabine2006}
S.~Zaglmayr.
\newblock {\em High Order Finite Element Methods for Electromagnetic Field
  Computation.}
\newblock PhD thesis, Graz University of Technology, 2006.

\bibitem{Zhang2019H}
Q.~Zhang, L.~Wang, and Z.~Zhang.
\newblock H($\text{curl}^2$)-conforming finite elements in 2 dimensions and
  applications to the quad-curl problem.
\newblock {\em SIAM Journal on Scientific Computing}, 41(3):A1527--A1547, 2019.

\bibitem{Zhang2018M2NA}
S.~Zhang.
\newblock Mixed schemes for quad-curl equations.
\newblock {\em Mathematical Modelling and Numerical Analysis}, 52(1):147--161,
  2018.

\bibitem{Zhang2017Regular}
S.~Zhang.
\newblock Regular decomposition and a framework of order reduced methods for
  fourth order problems.
\newblock {\em Numerische Mathematik}, 138(1):241--271, 2018.

\bibitem{Zheng2011A}
B.~Zheng and J.~Xu.
\newblock A nonconforming finite element method for fourth order curl equations
  in $\mathbb{R}^3$.
\newblock {\em Mathematics of Computation}, 80(276):1871--1886, 2011.

\end{thebibliography}
\end{document}